\documentclass{article}
\usepackage{epsfig}
\usepackage{amsmath, amssymb, amsthm}	
\usepackage{color}			
\usepackage{graphicx}
\usepackage{caption}
\usepackage{subcaption}

\newcommand{\be}{\begin{equation}}
\newcommand{\ee}{\end{equation}}
\newcommand{\ba}{\begin{eqnarray}}
\newcommand{\ea}{\end{eqnarray}}

\newcommand{\la}{\label}
\newcommand{\ep}{\Delta t}

\def\X{{\bf X}}

\newcommand{\R}{{I\!\!R}}
\newcommand{\N}{{I\!\!N}}

\def\R{{\rm I}\! {\rm R}}

\def\X{{\bf X}}

\newtheorem{remark}{Remark}[section]
\newtheorem{example}{Example}[section]

\begin{document}

\pagestyle{headings}

\title{Effective Simulation Methods for Structures with Local Nonlinearity: Magnus integrator and Successive Approximations}
\author{J\"urgen Geiser\thanks{Ruhr University of Bochum, Department of Electrical Engineering and Information Technology, Universit\"ats­str. 150, D-44801 Bochum, Germany, E-mail: juergen.geiser@ruhr-uni-bochum.de}
  \and Vahid Yaghoubi \thanks{Chalmers University of Technology, Department of Applied Mechanics, Gothenburg, Sweden, E-mail: yaghoubi@chalmers.se}}

\maketitle

\begin{abstract}

In the following, we discuss nonlinear simulations of nonlinear dynamical
systems, which are applied in technical and biological models.
We deal with different ideas to overcome the treatment of the nonlinearities
and discuss a novel splitting approach.
While Magnus expansion has been intensely studied and 
widely applied for solving explicitly time-dependent problems,
it can also be extended to nonlinear problems. By the way it is 
delicate to extend, while an exponential character
have to be computed.
Alternative methods, like successive approximation
methods, might be an attractive tool, which take into account 
the temporally in-homogeneous equation (method of Tanabe and Sobolevski).
In this work, we consider nonlinear stability analysis with 
numerical experiments and compare standard integrators to
our novel approaches.
\end{abstract}



\section{Introduction}

In this paper, we concentrate on solving
nonlinear evolution equations, which arose in nonlinear
dynamical applications, e.g., biological growth regimes or plasma
simulations.
We apply the following nonlinear differential equation,
\begin{eqnarray}
\label{equ1}
&& \partial_t \; u = A_1 u + A_2(u), \; u(0) = u_0 ,
\end{eqnarray}
where $A_1$ is the linear operator, while $A_2(u)$ is the nonlinear
operator.

To solve such delicate nonlinear differential equations, we have different
approaches:
\begin{itemize}
\item Approximation of the nonlinear term via time-dependent terms (Magnus expension).
\item Approximation of the nonlinear terms via multiscale expensions (successive approximation).
\end{itemize}
The Magnus expansion \cite{blan08} is an attractive and widely applied
method of solving explicitly time-dependent problems.
However, it requires computing
time-integrals and nested commutators to higher orders.
Successive approximation is based on recursive integral formulations
in which an iterative method is enforce the time dependency.

The paper is outlined as follows:
In Section \ref{intro}, we summarizes the Magnus expansion
and its application to Hamiltonian systems. 
Further, we show how AB-, Verlet and successive approximation method can be applied to
any exponential-splitting algorithms in Section \ref{nonlinear_1}.
In Section \ref{apply} we present the numerical results of the
splitting schemes.
In Section \ref{concl}, we briefly summarize our results.

\section{Introduction to splitting methods}
\label{intro}

For nonlinear problems, the applications of splitting methods are more
delicate because of resolving the nonlinear operators.
We apply a nonlinear approach based on Magnus expansion and successive approximations.

We concentrate on approximation to the solution of the 
nonlinear evolution equation, e.g. 
time-dependent Schr\"odinger equation,
\begin{eqnarray}
\label{equ1_1}
&&  \partial_t \; u = F(u(t)) + g(u) = A u(t) + B(u(t)) u(t) + g(u),  \; u(0) = u_0 ,
\end{eqnarray}
with the unbounded operators $A: D(A) \subset \X \rightarrow \X$
and  $B: D(B) \subset \X \rightarrow \X$.
We have further $F(v) = A(v) + B(v), v \in D(A) \cap D(B) $
and $g(u)$ is a nonlinear function.

We assume to have suitable chosen sub-spaces of the underlying
Banach space $(X, || \cdot ||_{X})$ such that $ D(F) = D(A) \cap D(B) \neq \emptyset  $

The exact solution of the evolution problem \ref{equ1_1} is given as:
\begin{eqnarray}
\label{equ1}
&&  u(t) = {\cal E}_F(t, u(0)) + \int_0^t {\cal E}_F(t - \tau, g(u(\tau)) d\tau , \; 0 \le t \le T ,
\end{eqnarray}
with the evolution operator ${\cal E}_F$ depending on the actual time $t$
and the initial value $u(0)$.

\begin{example}
For the linear case, means $F(u) = A u + B u$ the 
evolution operator is given as:
\begin{eqnarray}
\label{equ1}
&&  u(t) = \exp( t D_F)  u(0) + \int_0^t \exp((t-\tau) D_F) g(u(\tau)) d\tau , \; 0 \le t \le T ,
\end{eqnarray}
with $D_F = A + B$.

\end{example}

In the next subsections, we introduce the underlying splitting methods.

\section{Nonlinear Splitting Method}
\label{nonlinear_1}

We apply the abstract standard splitting schemes to the
multiproduct decomposition.

We have to carry out the following steps:
\begin{itemize}
\item Apply the nonlinear Strang splitting scheme,
\item embed the Strang splitting scheme into the
multiproduct expansion.
\end{itemize}

To apply the abstract setting of a nonlinear Magnus expansion, we deal with the following modified nonlinear equation
\be
\partial_t u = A u(t) + B(t, u(t)) u(t), \qquad u(0) = u_0 ,
\label{equ1_non}
\ee
where $A, B(t, u) \in [0, T] \times \X$ are non-commuting operators,
$\X$ is a general Banach space, e.g., $\X \subset \R^m$, where $m$ is the rank of the matrices.

To apply the nonlinear Magnus expansion, we deal with:
\be
\partial_t u = B(t, u(t)) u(t), \qquad u(0) = u_0 , \\
 u(t) = \exp(\Omega_{B}(t, u_0)) u_0 , \\
\la{equ1}
\ee
where the first order Magnus operator is given by
Euler's formula:
\be
\label{mag_1}
\Omega_{B,1}(t, u_0) = \int_{0}^{t} B(s, u(s)) ds = t \; B(0, u_0), 
\ee
and the second order Magnus operator is 
given by the midpoint rule:
\begin{eqnarray}
\label{mag_2}
\Omega_{B,2}(t, u_0) & = & \int_{0}^{t} B(s, \exp(\Omega_{B,1}(s, u_0)) ) ds, \\
 & = & t B(\frac{t}{2}, \exp( \Omega_{B,1}(\frac{t}{2}, u_0) ) )  \nonumber \\
 & = & t B(\frac{t}{2}, \frac{t}{2} B(0, u_0)),  \nonumber
\end{eqnarray}

or Trapezoidal-rule:
\begin{eqnarray}
\label{mag_2}
\Omega_{B,2}(t, u_0) & = & \int_{0}^{t} B(s, \exp(\Omega_{B,1}(s, u_0)) ) ds, \\
 & = & t/2  \left( B(0, u_0) + B(t, \exp( t B(0, u_0) u_0) ) \right) . \nonumber 
\end{eqnarray}

We can generalize the schemes with respect to more additional terms
to higher order schemes, see \cite{casas2006}

\begin{itemize}

 \item We apply the following A-B splitting scheme:

\begin{eqnarray}
&& \partial_t u_1 = A u_1(t) ,  \qquad u(t^n) = u_n , \\
&& \partial_t u_2 = B(t, u_2(t)) u_2(t), \qquad u_2(t^n) = u_1(t^{n+1}) ,
\label{equ1_non}
\end{eqnarray}
where the time-step is $\Delta t= t^{n+1} - t^n$ and the next 
solution is: \\
 $u(t^{n+1}) = u_2(t^{n+1})$.

Here we apply the kernels:
\begin{eqnarray}
&& u_2(t^{n+1}) = \exp(\Omega_{B,2}(\Delta t, u_1(t^{n+1})))  u_1(t^{n+1}), \\
&& u_1(t^{n+1}) = \exp(A \Delta t) u(t^n)
\end{eqnarray}

We have:
\begin{eqnarray}
{\cal T}_{1, BA}(\ep) u(t^n) & = &  \exp(\Omega_{B,2}(\Delta t, \exp(A \Delta t) u(t^n))) \exp(A \Delta t) u(t^n) , 
\end{eqnarray}

\item Verlet Splitting (Strang-Splitting)

We apply two A-B splitting, means:
\begin{eqnarray}
{\cal T}_{1, AB}(\ep) u(t^n) & = & \exp(A \Delta t)  \exp(\Omega_{B,1}(\Delta t, u(t^n))) , \\
{\cal T}_{1, BA}(\ep) u(t^n) & = &  \exp(\Omega_{B,2}(\Delta t, \exp(A \Delta t) u(t^n))) \exp(A \Delta t) u(t^n) ,
\end{eqnarray}
where $\Omega_{A, i}, \Omega_{B, i}$ are the Magnus expansions, see equations (\ref{mag_1}) and (\ref{mag_2}), of order $i$
and we obtain:
the scheme based on the symmetrical splitting by
\begin{eqnarray}
\label{linearized}
{\cal T}_{2}(\ep)u(t^n) & =& {\cal T}_{1,B A}(\frac{\ep}{2}) {\cal T}_{1, AB}(\frac{\ep}{2})u(t^n) . 
\end{eqnarray}

\item Standard Successive Approximation via linear operator $A$ (without multiscale approximation)

We deal with the equation:
\begin{eqnarray}
 && \frac{\partial u(t)}{\partial t} = A u(t) \; + \; B(t, u(t)) u(t), \; \mbox{with} \; \; u(0) = u_{init}(0) . 
\end{eqnarray}

Then the successive equations are given as:
\begin{eqnarray}
\label{mm_1}
&&  \frac{\partial u_0(t)}{\partial t} = A u_0(t) , \; u_0(t^n) = u(t^n) ,  \; t \in [t^n, t^{n+1}] , \\
&&\hspace{-1cm}  \frac{\partial u_1(t)}{\partial t} = A u_1(t) + \; B(u_0(t)) u_0(t)   ,  \; u_1(t^n) = u(t^n) ,  \; t \in [t^n, t^{n+1}]  ,\\  
&&  \hspace{-1cm}  \frac{\partial u_2(t)}{\partial t}  =  A u_2(t) + \; B(u_2(t)) u_1(t) ,  \; u_0(t^n) = u(t^n) , \; t \in [t^n, t^{n+1}] ,  \\
&& \hspace{-1cm}   \frac{\partial u_3(t)}{\partial t}  =  A u_3(t)  + \; B(u_2(t)) u_2(t) ,  \; u_3(t^n) = u(t^n) ,  \; t \in [t^n, t^{n+1}]  , \\
&& \vdots \nonumber
\end{eqnarray}
where $\Delta t = t^{n+1} - t^n$, we start with $u(0) = u_{init}(0)$,
the successive solution at $t^{n+1}$ with $i=3$ steps are given as
$u(t^{n+1}) = u_3(t^{n+1})$.

By integration, we have the following solutions, we start with $n=0$ and $u(0) = u(t^0) = u_{init}(0)$:
\begin{eqnarray}
\label{mm_1}
  u_0(t^{n+1}) & = & \exp(A \Delta t) u(t^n) , \\
  u_1(t^{n+1}) & = & \exp(A \Delta t) u(t^n) \nonumber \\
             &&  + \int_{t^n}^{t^{n+1}} \exp(A (t^{n+1}- s)) \; B(u_0(s)) \; u_0(s) \; ds   , \\
  u_2(t^{n+1}) & = & \exp(A \Delta t) u(t^n) \nonumber \\
             && + \int_{t^n}^{t^{n+1}} \exp(A (t^{n+1} - s)) B(u_1(s)) u_1(s) \; ds  ,
\end{eqnarray}
with $\Delta t = t^{n+1} - t^n$ is the time step.

We apply a simple trapezoidal-rule to the integrals and obtain:
\begin{eqnarray}
\label{mm_1}
&&  \tilde{u}_1(t^{n+1}) = \frac{\Delta t}{2}  \left( B(u_0(t^{n+1})) \; u_0(t^{n+1}) \right. \\
&& \left. + \exp(A \Delta t) B(u_0(t^n)) \; u_0(t^n) \right)   ,\nonumber \\
&&  \tilde{u}_2(t^{n+1}) = \frac{\Delta t}{2} \left( B(u_1(t^{n+1})) u_1(t^{n+1}) \right. \\
&&\left. +  \exp(A \;  \Delta t) B(u_1(t^{n})) u_1(t^{n}) \right)  , \nonumber
\end{eqnarray}
where we applied mid-point rule or Simpson's-rule,
then we have for $3$ iterative steps, $u(t^{n+1})$ is given as
\begin{eqnarray}
u(t^{n+1}) = u_0(t^{n+1}) + \tilde{u}_1(t^{n+1}) + \tilde{u}_2(t^{n+1}) .
\end{eqnarray}

\item Successive Approximation via linear operator $A$ (Multiscale expansion)

We deal with the multiscale idea:
\begin{eqnarray}
 && \frac{\partial u(t)}{\partial t} = A u(t) \; + \; \epsilon \; B(t, u(t)) u(t), \; \mbox{with} \; \; u(0) = u_{init}(0) , 
\end{eqnarray}
where $0 < \epsilon \le 1$. For $\epsilon \rightarrow 1$, we have the original equation.

We derive a solutions $u(t, \epsilon)$ and apply:
\begin{eqnarray}
u(t, \epsilon) = u_0(t) + \epsilon u_1(t) + \epsilon^2 u_2(t) + \ldots +\epsilon^J u_J(t)  ,
\end{eqnarray}
with the initial conditions  $u(0, \epsilon) = u(0)$ and $J \in \N^+$ is a fixed iteration number.

Then the hierarchical equations are given as:
\begin{eqnarray}
\label{mm_1}
  \frac{\partial u_0(t)}{\partial t} & = & A u_0(t) , \; u_0(t^n) = u(t^n) , t \in [t^n, t^{n+1}] , \\
  \frac{\partial u_1(t)}{\partial t} & = & A u_1(t) + \; B(u_0(t)) u_0(t)   , \\
&&  \; u_1(t^n) = u(t^n) , t \in [t^n, t^{n+1}] , \nonumber  \\  
  \frac{\partial u_2(t)}{\partial t} & = & A u_2(t) \\
 && + \; B(u_0(t)) u_1(t) + \; B'(u_0(t)) u_1(t) u_0(t) , \nonumber \\
&& \; u_2(t^n) = u(t^n) , t \in [t^n, t^{n+1}] , \nonumber \\
  \frac{\partial u_3(t)}{\partial t} & = & A u_3(t)  + \; B(u_0(t)) u_2(t) + \; B'(u_0(t)) u_2(t) u_0(t) \\
&& + \; B'(u_0(t)) u_1(t) u_1(t) + \; B''(u_0(t)) u_1(t) u_1(t) u_0(t) , \nonumber \\
&& \; u_3(t^n) = u(t^n) , t \in [t^n, t^{n+1}] , \nonumber \\
& \vdots & \nonumber
\end{eqnarray}
where we extend $B( u_0(t) + \epsilon (u_1(t) + \epsilon u_2(t) + \ldots)) = B(u_0(t)) + \epsilon (u_1(t) + \epsilon u_2(t) + \ldots) B'(u_0) + (\epsilon (u_1(t) + \epsilon u_2(t) + \ldots))^2 B''(u_0) * \ldots$ with $B'(u_0) = \frac{\partial B(u)}{\partial u}|_{u = u_0}$ and   $B''(u_0) = \frac{\partial^2 B(u)}{\partial u^2}|_{u = u_0}$.
We have also to expand the initial conditions to $u_0(0) = u_{init}(0)$ and $u_j(0) = 0, \; \forall j = 1, \ldots, J$.

By integration, we have the following solutions, we start with $n=0$ and $u(0) = u(t^0) = u_{init}(0)$:
\begin{eqnarray}
\label{mm_1}
  u_0(t^{n+1}) & = & \exp(A \Delta t) u(t^n) , \\
  u_1(t^{n+1}) & = & \exp(A \Delta t) u(t^n) + \int_{t^n}^{t^{n+1}} \exp(A (t^{n+1}- s)) \; B(u_0(s)) \; u_0(s) \; ds   , \\
  u_2(t^{n+1}) & = & \exp(A \Delta t) u(t^n) \\
&& + \int_{t^n}^{t^{n+1}} \exp(A (t^{n+1} - s)) \left( B(u_0(s)) u_1(s) + \; B'(u_0(s)) u_1(s) u_0(s) \right) \; ds  , \nonumber
\end{eqnarray}
with $\Delta t = t^{n+1} - t^n$ is the time step.

We apply a simple trapezoidal-rule to the integrals and obtain:
\begin{eqnarray}
\label{mm_1}
&&  \tilde{u}_1(t^{n+1}) = \frac{\Delta t}{2}  \left( B(u_0(t^{n+1})) \; u_0(t^{n+1}) + \exp(A \Delta t) B(u_0(t^n)) \; u_0(t^n) \right)   , \\
&&  \tilde{u}_2(t^{n+1}) = \frac{\Delta t}{2} \left( \left( B(u_0(t^{n+1})) u_1(t^{n+1}) + \; B'(u_0(t^{n+1})) u_1(t^{n+1}) u_0(t^{n+1}) \right) \right. \nonumber \\
&& \left. +  \exp(A \;  \Delta t) \left( B(u_0(t^{n})) u_1(t^{n}) + \; B'(u_0(t^{n})) u_1(t^{n}) u_0(t^{n}) \right) \right)  , 
\end{eqnarray}
where we also can apply mid-point rule or Simpson's-rule,
then we have for $u(t^{n+1})$
\begin{eqnarray}
u(t^{n+1}) = u_0(t^{n+1}) + \tilde{u}_1(t^{n+1}) + \tilde{u}_2(t^{n+1}) + \tilde{u}_3(t^{n+1})  ,
\end{eqnarray}

\item Successive Approximation via linear operator $B$ (multiscale expansion)

We deal with the multiscale idea:
\begin{eqnarray}
 && \frac{\partial u(t)}{\partial t} = \epsilon \; A u(t) \; + \; B(t, u(t)) u(t), \; \mbox{with} \; \; u(0) = u_{init}(0) , 
\end{eqnarray}
where $0 < \epsilon \le 1$. For $\epsilon \rightarrow 1$, we have the original equation.

We derive a solutions $u(t, \epsilon)$ and apply:
\begin{eqnarray}
u(t, \epsilon) = u_0(t) + \epsilon u_1(t) + \epsilon^2 u_2(t) + \ldots +\epsilon^J u_J(t)  ,
\end{eqnarray}
with the initial conditions  $u(0, \epsilon) = u(0)$ and $J \in \N^+$ is a fixed iteration number.

Then the hierarchical equations are given as:
\begin{eqnarray}
\label{mm_1}
  \frac{\partial u_0(t)}{\partial t} & = & B(u_0(t)) u_0(t) , \; u_0(t^n) = u(t^n) , \\
  \frac{\partial u_1(t)}{\partial t} & = &  \left( B(u_0(t)) + \; B'(u_0(t)) u_0(t) \right) u_1(t) +  A u_0(t)   , \\
&&  \; u_1(t^n) = u(t^n) , \nonumber \\  
  \frac{\partial u_2(t)}{\partial t} & = &   \; \left( B(u_0(t)) u_2(t) + \; B'(u_0(t)) u_0(t) \right) u_2(t) \\
&& + \; B'(u_0(t)) u_1(t) u_1(t) + \; B''(u_0(t)) u_1(t) u_1(t) u_0(t)  + A u_1(t) ,  \nonumber \\
&& \; u_2(t^n) = u(t^n) , \nonumber \\
& \vdots & \nonumber
\end{eqnarray}
where the initialization is $u_0(0) = u_{init}(0)$ and we extend $B( u_0(t) + \epsilon (u_1(t) + \epsilon u_2(t) + \ldots)) = B(u_0(t)) + \epsilon (u_1(t) + \epsilon u_2(t) + \ldots) B'(u_0) + (\epsilon (u_1(t) + \epsilon u_2(t) + \ldots))^2 B''(u_0) * \ldots$ with $B'(u_0) = \frac{\partial B(u)}{\partial u}|_{u = u_0}$ and \\   $B''(u_0) = \frac{\partial^2 B(u)}{\partial u^2}|_{u = u_0}$.
We have also to expand the initial conditions to $u_0(0) = u_{init}(0)$ and $u_j(0) = 0, \; \forall j = 1, \ldots, J$.

By integration, we have the following solutions, we start with $n=0$ and $u(0) = u(t^0) = u_{init}(0)$:
\begin{eqnarray}
\label{mm_1}
&&  u_0(t^{n+1}) = \exp(B(u(t^n)) \Delta t) u(t^n) , \\
&&  u_1(t^{n+1}) = \exp( B(u_0(t^n)) + B'(u_0(t^n)) \; u_0(t^n) )  u(t^n) \nonumber \\
&&  + \int_{t^n}^{t^{n+1}} \exp\left(\int_{s_1}^{t^{n+1}} B(u_0(s_1)) + B'(u_0(s_1)) \; u_0(s_1) \; ds_1 \right) \; A \; u_0(s) \; ds   , \nonumber \\
&& = \exp( B(u_0(t^n)) + B'(u_0(t^n)) \; u_0(t^n) )  u(t^n)  \\
&&  +  \int_{t^n}^{t^{n+1}} \exp \left( \left(B(u_0(s)) +  B'(u_0(s)) \; u_0(s) \right)\; (t^{n+1} - s) \right) \; A \; u_0(s) \; ds  \nonumber ,
\end{eqnarray}
with $\Delta t = t^{n+1} - t^n$ is the time step.

We apply a simple trapezoidal-rule to the integrals and obtain:
\begin{eqnarray}
\label{mm_1}
&&  u_1(t^{n+1}) = \frac{\Delta t}{2} \left( A u_0(t^{n+1}) \right. \nonumber \\
&& \left. +  \exp \left( \left(B(u_0(t^n)) +  B'(u_0(t^n)) \; u_0(t^n) \right)\; \Delta t \right) \; A \; u_0(t^{n})  \right) ,
\end{eqnarray}
where we also can apply mid-point rule or Simpson's-rule,
then we have for $u(t^{n+1})$
\begin{eqnarray}
u(t^{n+1}) = \tilde{u}_0(t^{n+1}) + \tilde{u}_1(t^{n+1}) = u_1(t^{n+1}) ,
\end{eqnarray}

\end{itemize}

\begin{remark}

The benefit of the A-B and Strang-Splitting schemes are based on the
explicit Magnus expansion and the fully decomposition of operator $A$ and $B$.

The benefit of the successive approximation scheme is the idea to skip the 
Magnus expansion and to deal with e relaxation over $B$. Here we have a weakly coupling based on the frozen solution (previous iterated solution)
and we could deal only with matrix multiplications.
\end{remark}

\begin{remark}
While AB splitting has the idea to decompose into an A and B operator-equation,
we have to compute for the A-equation $\exp(A)$ and for the B-equation 
the nonlinear $\exp(\int_{t^n}^{t^{n+1}} B(u(t)) \; dt)$ part.
The last term is expensive.
The successive approximation has the following idea to iterate or perturb
only to the A-operator, means we have the cheap part $\exp(A)$ and the
integral part is only with the right-hand $B(u_0(t))$, that is given of the
previous solution and there is no need to apply the Magnus-integrators
for the exponential expension.
\end{remark}

\section{Numerical Experiments}
\label{apply}

In the following section, we deal with experiments to
verify the benefit of our methods.
At the beginning, we propose introductory examples to
compare the methods. In the next examples, applications to
nonlinear differential equations, as Bernoulli's equation
and Fisher's equation for biological models.

\subsection{First test example of a nonlinear ODE: Bernoulli's equation}

We deal with a nonlinear ODE (Bernoulli's equation)
and split it into linear and nonlinear operators.

First we examine the non linear Bernoulli-Equation\\
\begin{eqnarray}
\label{num7}
\frac{\partial u(t)}{\partial t} &=& \lambda_1 u(t) + \lambda_2 u^n(t) \; , \\
u(0) &=& 1 \; , 
\end{eqnarray}
with the analytic solution
\begin{eqnarray}
\label{analyt}
 u(t) = \left[ (1+\frac{\lambda_2}{\lambda_1})\exp(\lambda_1t(1-n))-\frac{\lambda_2}{\lambda_1})\right]^{\frac{1}{1-n}} \; .
\end{eqnarray}

For the computations, we choose $n=2$ , $\lambda_1 = -1$, $\lambda_2 = -1, -2, -10$ and $\Delta t=10^{-2}$ . \\

We rewrite the equation-system (\ref{num7}) in
operator notation, and obtain the following equations :
\begin{eqnarray}
\label{num_8}
  \partial_t u & = & A u + B(u) u \; , \\
   u(0) & = & u_0 ,
\end{eqnarray}
where $ u(t)=(u_1(t),u_2(t))^T$ for $t \in [0,T]$.

Our split operators are
\begin{equation}
\label{num_9}
A u = \lambda_1 u \; , \\
B(u) u = \lambda_2 \; u^{m-1} u\; ,
\end{equation}
with $m = 2$.

We also have a non-commutative behavior of the
nonlinear operators, means $[A(u),B(u)] = \left( \frac{\partial B(u)}{\partial u} A u - \frac{\partial A u}{\partial u} B(u) \right) = \lambda_2 (m-1) u^{m-2} \lambda_2 u - \lambda_1 \lambda_2 u^{m-1} \neq 0$.

We have the following results with the $L_2$ and $L_{inf}$-error of our scheme, related to the analytic solution (\ref{analyt}) in Figure \ref{error_1} and \ref{error_2}.

\begin{table}[h]
\begin{center}
\begin{tabular}{||c|c|c|c||}
\hline \hline
Numerical Method & $err_{L_2}$ & $err_{L_{inf}}$ & Comput. time [sec] \\
\hline
AB-Splitting & $0.962$  & $0.183$ & $0.030$ \\
\hline
Strang-splitting & $0.752$ & $0.118$ & $0.055$ \\
\hline
Standard Successive & $0.283$ & $0.100$ & $0.013$ \\
\hline
Multiscale Successive & $0.098$ & $0.035$ & $0.027$ \\
\hline 
\end{tabular}
\caption{\label{error_1} Numerical error of the different numerical methods for the Bernoulli's equation.}
\end{center}
\end{table}

\begin{table}[h]
\begin{center}
\begin{tabular}{||c|c|c|c|c||}
\hline \hline
Numerical Method & $\Delta t$& $err_{L_2}$ & $err_{L_{inf}}$ & Comput. time [sec] \\
\hline
Stand. Successive & $0.02$ & 0.260 & 0.096 & 0.002 \\
Multi. Successive & $0.02$ & 0.087 & 0.032 & 0.004 \\
AB-Splitting & $0.02$ &   0.800 & 0.173 & 0.017 \\
Strang-splitting & $0.02$ & 0.588 & 0.110 & 0.011 \\
\hline
Stand. Successive & $0.01$ & 0.160 & 0.039 & 0.004 \\
Multi. Successive & $0.01$ & 0.026  & 0.007 & 0.009 \\
AB-Splitting      & $0.01$ & 1.154  & 0.177 & 0.011 \\ 
Strang-splitting  & $0.01$ & 0.826  & 0.109 & 0.018 \\
\hline
Stand. Successive & $0.005$ & 0.106 & 0.018 & 0.008 \\
Multi. Successive & $0.005$ & 0.008 & 0.002 & 0.016 \\
AB-Splitting      & $0.005$ & 1.648 & 0.179 & 0.020 \\
Strang-splitting  & $0.005$ & 1.164 & 0.109 & 0.037 \\
\hline
Stand. Successive & $0.0025$ & 0.073 & 0.009 & 0.016 \\
Multi. Successive & $0.0025$ & 0.003 & 0.000 & 0.034 \\
AB-Splitting      & $0.0025$ & 2.342 & 0.180 & 0.042 \\
Strang-splitting  & $0.0025$ & 1.642 & 0.109 & 0.075 \\
\hline 
\end{tabular}
\caption{\label{error_2} Numerical error of the different numerical methods for the Bernoulli's equation.}
\end{center}
\end{table}

We apply the convergence-rates as
\begin{eqnarray}
\label{kap7_gleich3}
\rho_{L_{2,\Delta t}}(t) =
\frac{log \bigg(\frac{E_{L_{2,\Delta t/2}}(t)}{E_{L_{2,\Delta t}}(t)}\bigg)}{log(0.5)} .
\end{eqnarray}

We obtain a stagnation of the numerical errors of the AB- and Strang-splitting,
means, we could not obtain a convergent behavior. Instead with the standard and multiscale successive approach, we obtain a convergence with following convergence rates.

The convergence-rates of the different schemes are given in Table \ref{table_2}:
\begin{table}[ht]
\begin{center}
\begin{tabular}{|c || c | c | }
\hline
time-step & Standard Succ. $(\rho_{L_2})$ &  Multi-Splitt $(\rho_{L_2})$ \\
\hline
\hline 
$\Delta t$ & 0.700 & 1.7425 \\ 
$\Delta t/2$ & 0.594 & 1.70043 \\ 
$\Delta t/4$ & 0.5381 & 1.42 \\
\hline
\end{tabular}
\end{center}
\vspace{0.1cm}
\caption{\label{table_2} Order of convergence of $L_2$-error for the different numerical splitting schemes.}
\end{table}

\begin{figure}[ht]
\begin{center}  
\includegraphics[width=7.0cm,angle=-0]{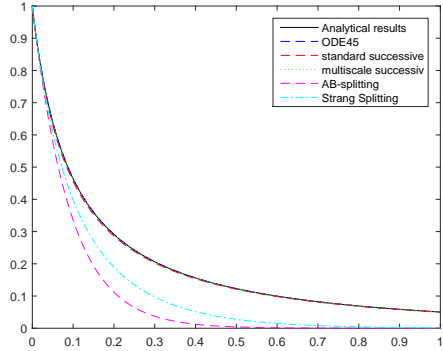} 
\end{center}
\caption{\label{figure_1_1} The solutions of the Bernoulli's equation solved with the different splitting schemes.}
\end{figure}

The results of the different schemes are shown in Figure \ref{figure_1_1}.

%

\begin{remark}
We compare the different standard splitting scheme, e.g., A-B splitting, Strang-splitting and the standard successive approximation based on a Bernoulli's equation with a moderate nonlinearity. We obtain the best results in the case of the novel multiscale successive approximation method, while we include a more accurate resolution of the fast scales. We also obtain an fast method, that is competitive with the simple AB-splitting scheme. Such we improve the results with a more adapted novel scheme.
\end{remark}

\subsection{Second test example: Diffusion-Reaction equation with
nonlinear reaction (Fisher's equation)}

We deal with the Fisher Equation, which describe the
spreading of genes see \cite{fisher1937} and has found applications in different
fields of research ranging from ecology \cite{kolomgorov1991} to plasma
physics \cite{gilding2004}.

We deal with a nonlinear PDE and split it into linear and nonlinear operators,
while we can compare to a analytical solution.

The Fisher's equation is given as
\begin{eqnarray}
\frac{\partial u(x, t)}{\partial t} & = & D \partial_{xx} u + r (1 - \frac{u}{K}) u \; , \; \mbox{in} \; (x, t) \in \Omega \times [0, T] \\
u(x, 0) &=& g(x)  \; \mbox{on} \; x \in \Omega_{init} , \\
u(x, t) &=& u_{analyt}(x,t) \;   , \; \mbox{on} \; (x, t) \in \partial \Omega \times [0, T],
\end{eqnarray}
where we assume $\Omega \subset \R$.
$u$ is the solution function, the initial condition is $g(x)$.
Form the dynamical view-point, we apply a homogeneous medium with $D$ as diffusion coefficient and we embed a growth of a logistic function, see \cite{vander2010}
with the $r$ is the growth rate and $K$ is the carrying capacity.

The analytical solution is given as:
\begin{eqnarray}
\label{fisher_1}
 u_{analyt}(x, t) = \frac{\exp(r \; t)}{1 + \frac{1}{K} \left( \exp(r \; t) \tilde{u}(x, t) - g(x) \right)} \tilde{u}(x, t) \; .
\end{eqnarray}
where $\tilde{u}(x,t) = \frac{1}{\sqrt{1 + 4 D \; t}} \int_{-\infty}^{\infty} \exp(- \frac{(x - \sigma)^2}{1 + 4 D \; t})  \; g(\sigma) \; d\sigma$.

and we apply the case $g(x) = \exp(-x^2)$ and we have \\
$\tilde{u}(x,t) = \frac{1}{\sqrt{1 + 4 D \; t}} \exp(- \frac{x^2}{1 + 4 D \; t})$.

We rewrite the equation-system (\ref{fisher_1}) in
operator notation, and obtain the following equations :
\begin{eqnarray}
\label{num_8}
  \partial_t u & = & A u + B(u) u \; , \\
   u(0) & = & u_0 ,
\end{eqnarray}
and we split our operators to a linear and nonlinear one:
\begin{eqnarray}
\label{num_9}
A =   D \partial_{xx} + r \\
B(u) = - r \frac{u}{K}  ,
\end{eqnarray}
with $D = 0.01, \; r = 1, K = 1$, later we apply the multiscale case with $K = 0.5, 0.25$.

\begin{itemize}
\item Analytical Solution of the Diffusion-Reaction Part: \\
Based on the one-dimensional problem, we can apply the analytical solution of the diffusion-convection part means:
\begin{eqnarray}
\frac{\partial u(x, t)}{\partial t} & = & D \partial_{xx} u + r  u \; , \; \mbox{in} \; (x, t) \in \Omega \times [0, T] \\
u(x, 0) &=& u_0(x)  \; \mbox{on} \; x \in \Omega_{init} , \\
u(x, t) &=& u_{analyt}(x,t) \;   , \; \mbox{on} \; (x, t) \in \partial \Omega \times [0, T],
\end{eqnarray}
where we assume $\Omega \subset \R$.
$u$ is the solution function, the initial condition is $u_0(x)$.
Further $D$ is the diffusion coefficient and $r$ the growth rate.

The analytical solution is given as:
\begin{eqnarray}
\label{fisher_1}
 u_{analyt}(x, t) = \exp(r \; t) \tilde{u}(x, t) \; .
\end{eqnarray}
where $\tilde{u}(x,t) = \frac{1}{\sqrt{1 + 4 D \; t}} \int_{-\infty}^{\infty} \exp(- \frac{(x - \sigma)^2}{1 + 4 D \; t})  \; u_0(\sigma) \; d\sigma$.

\item Numerical Solution of the Diffusion part: \\

The operator $A$ is discretized as:
\begin{eqnarray}
A & = &  \frac{D}{\Delta x^2}\cdot  \left(\begin{array}{rrrrr}
 -2 & 1 & ~ & ~ & ~ \\
 1 & -2 & 1 & ~ & ~ \\
 ~ & \ddots & \ddots & \ddots & ~ \\
 ~ & ~ & 1 & -2 & 1 \\
 ~ & ~ & ~ & 1 & -2
\end{array}\right) \nonumber \\
& + & r  \left(\begin{array}{rrrrr}
 1 & 0 & ~ & ~ & ~ \\
 0 & 1 & ~ & ~ & ~ \\
 ~ & \ddots & \ddots & ~ & ~ \\
 ~ & ~ & 0 & 1 & 0 \\
 ~ & ~ & ~ & 0 & 1
\end{array}\right)~\in~\R^{I \times I},
\end{eqnarray}
where $I$ are the number of discretization points, e.g. $I = 100$.

The operator $B$ is written in the operator notation as:
\begin{eqnarray}
B(u) =  - \frac{r}{K}  \left(\begin{array}{rrrrr}
 u_1(t) & 0 & ~ & ~ & ~ \\
 0 & u_2(t) & ~ & ~ & ~ \\
 ~ & \ddots & \ddots & ~ & ~ \\
 ~ & ~ & 0 & u_{I-1}(t) & 0 \\
 ~ & ~ & ~ & 0 & u_I(t)
\end{array}\right)~\in~\R^{I \times I},
\end{eqnarray}
where $I$ are the number of discretization points, e.g. $I = 100$.

The solution is given as ${\bf u}(t) = (u_1(t), \ldots, u_I(t))$,
where $u_i(t) = u(x_i,t), x_i = \Delta x \; i, \; i = 1, \ldots, I$.

Now, we can apply the discretized operator equations
\begin{eqnarray}
\frac{\partial {\bf u}(t)}{\partial t} & = & A {\bf u} + B({\bf u}) {\bf u} \; , \; \mbox{in} \; t \in [0, T] \\
{\bf u}(0) &=& (g(x_1), \ldots, g(x_I))^t ,
\end{eqnarray}
 to our schemes.

\end{itemize}

The nonlinear behaves in a non-commutative
manner $A' B - B' A \neq 0$.

We apply finite differences or finite elements to the spatial operator.
With $\Omega = [0, 10]$, $\Delta x = 0.05$ \\

\begin{figure}[ht]
\begin{center}  
\includegraphics[width=5.0cm,angle=-0]{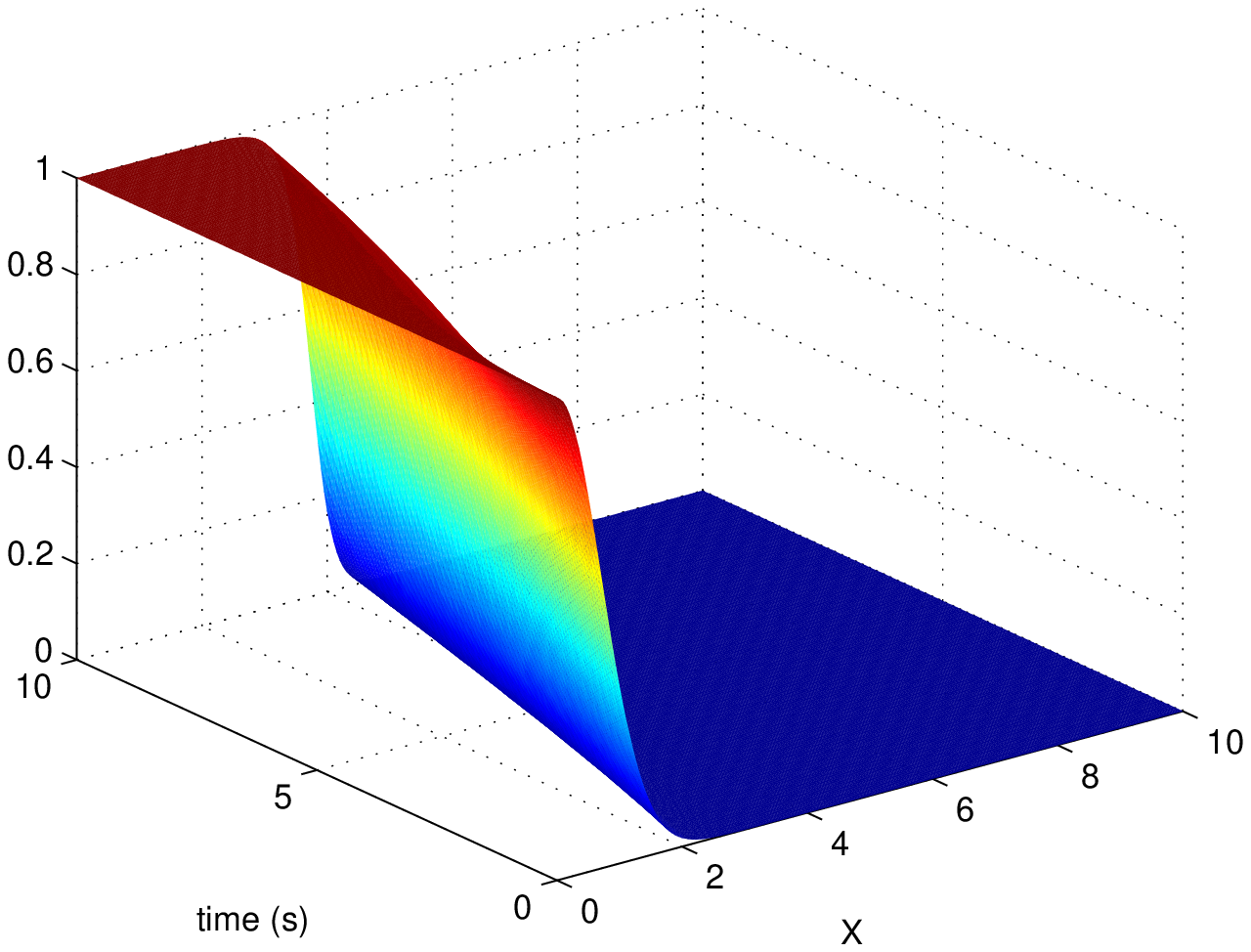} 
\includegraphics[width=5.0cm,angle=-0]{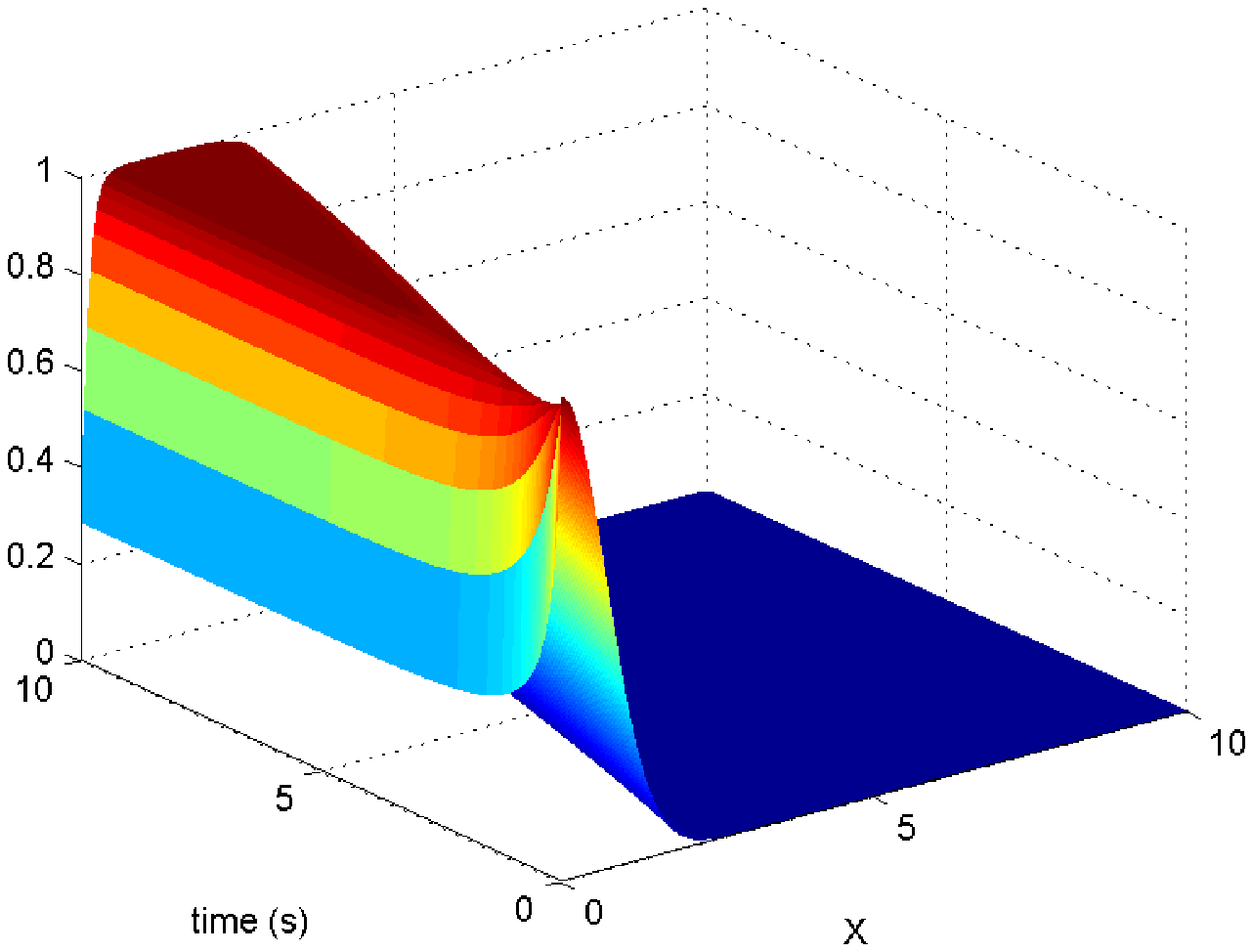} 
\end{center}
\caption{\label{figure_1} The solutions of the 1D Fisher's equation with the analytical and numerical splitting schemes.}
\end{figure}

We have the following results:

The results of the different schemes are shown in Figure \ref{figure_1}.

%

We have the following error in an $L_2$ Banach space:
\begin{eqnarray}
\label{kap7_gleich3}
E_{L_{2, \Delta x}}(t) & = & \sqrt{\int_{\Omega_h} (u_{ana}(x, t) - u_{num}(x, t))^2 \; dx } \\
              & = & \sqrt{\Delta x\sum_{i=1}^N (u_{ana}(x_i, t) - u_{num}(x_i, t))^2  } 
\end{eqnarray}

\begin{figure}[ht]
\begin{center}  
\includegraphics[width=3.5cm,angle=-0]{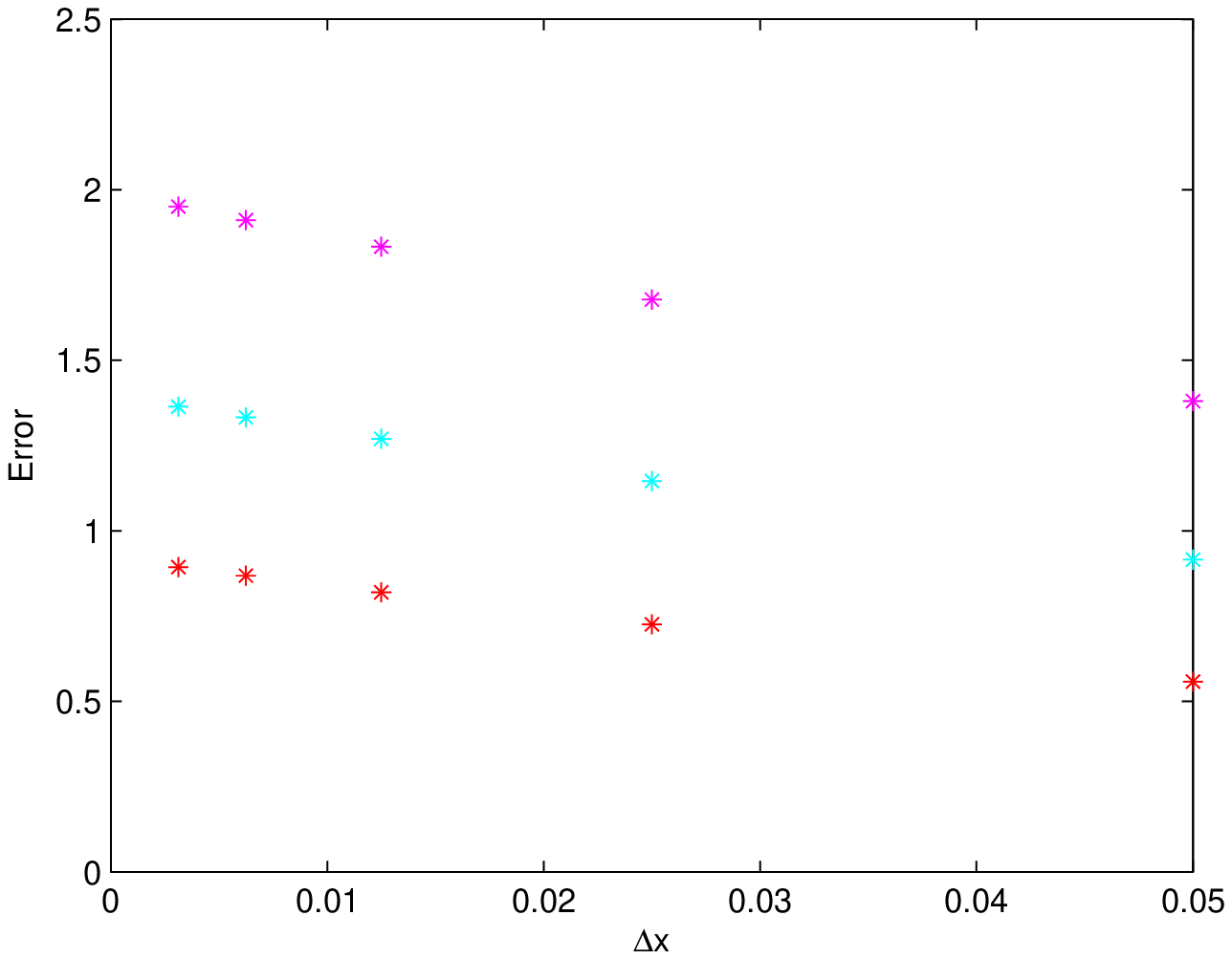} 
\includegraphics[width=3.5cm,angle=-0]{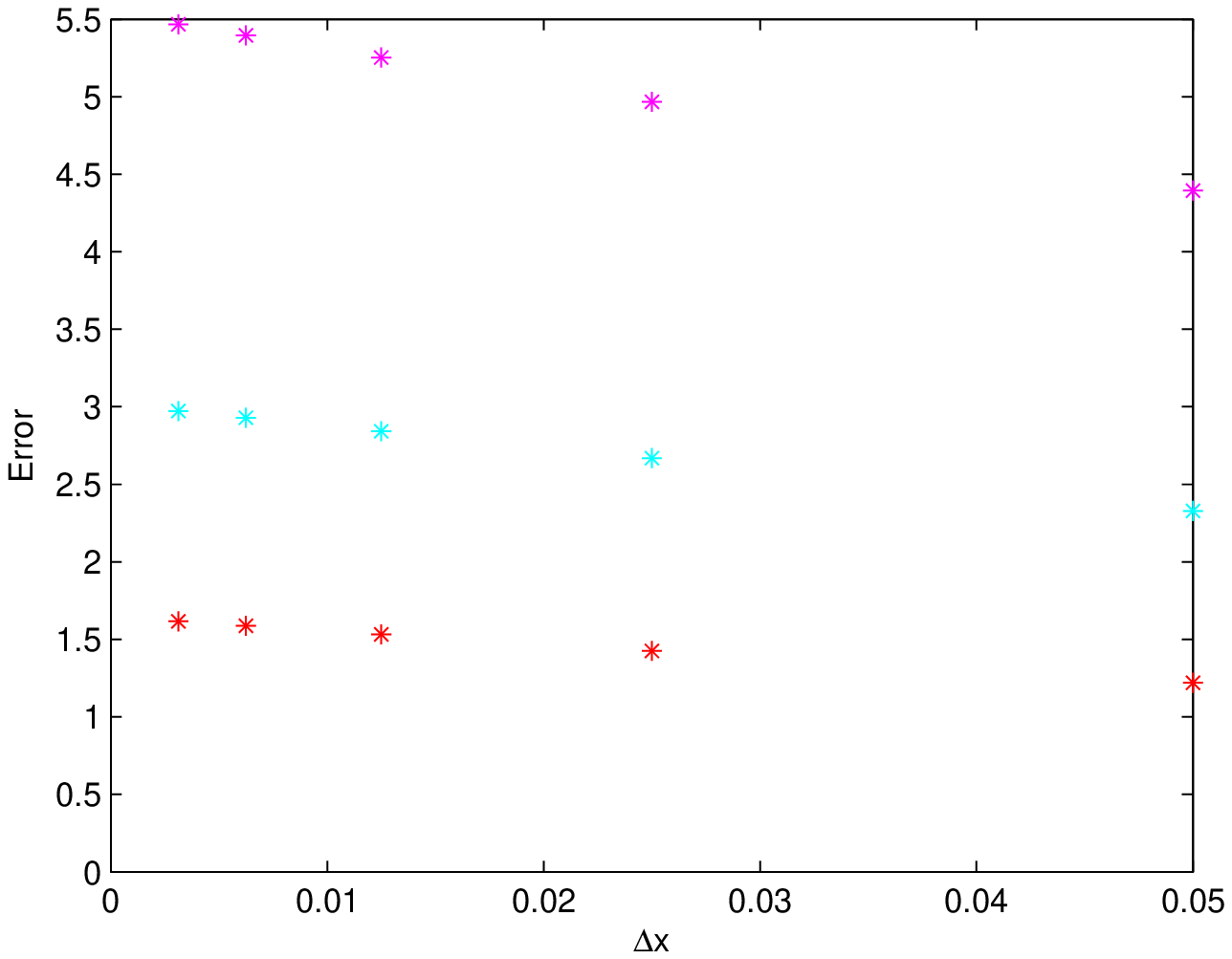} 
\includegraphics[width=3.5cm,angle=-0]{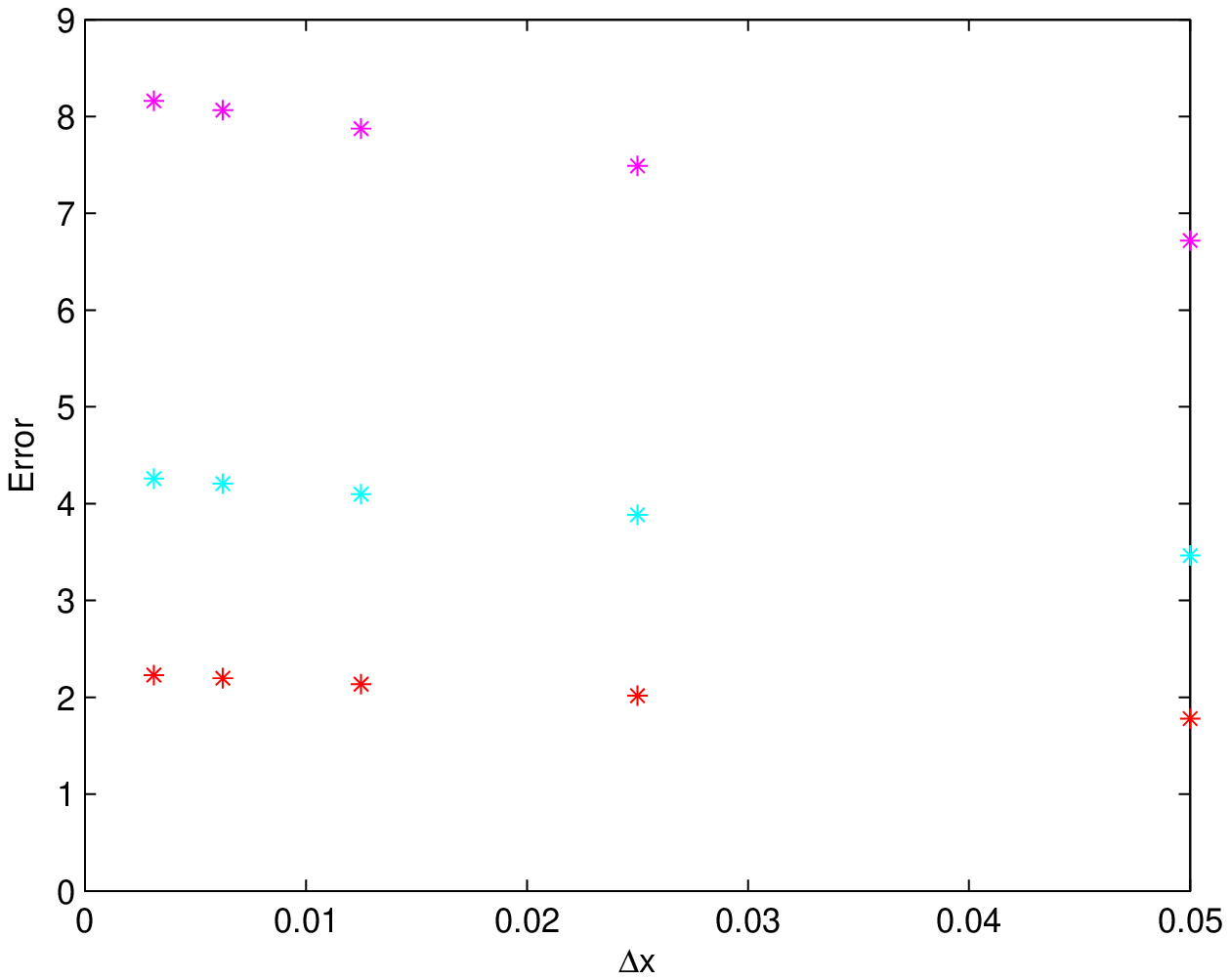} 
\end{center}
\caption{\label{figure_3} The solutions of the 1D Fisher's equation with the numerical splitting schemes for $K = 1, \; 0.5, \; \mbox{and} \; 0.25$ (asterisks in magenta $K=1.0$, cyan $K=0.5$ and red $K = 0.25$),  (left figure: $T=1$, middle figure: $T=5$ and right figure: $T=10$).}
\end{figure}

The results of the different schemes are shown in Figure \ref{figure_3}.

\begin{figure}[ht]
\begin{center}  
\includegraphics[width=3.5cm,angle=-0]{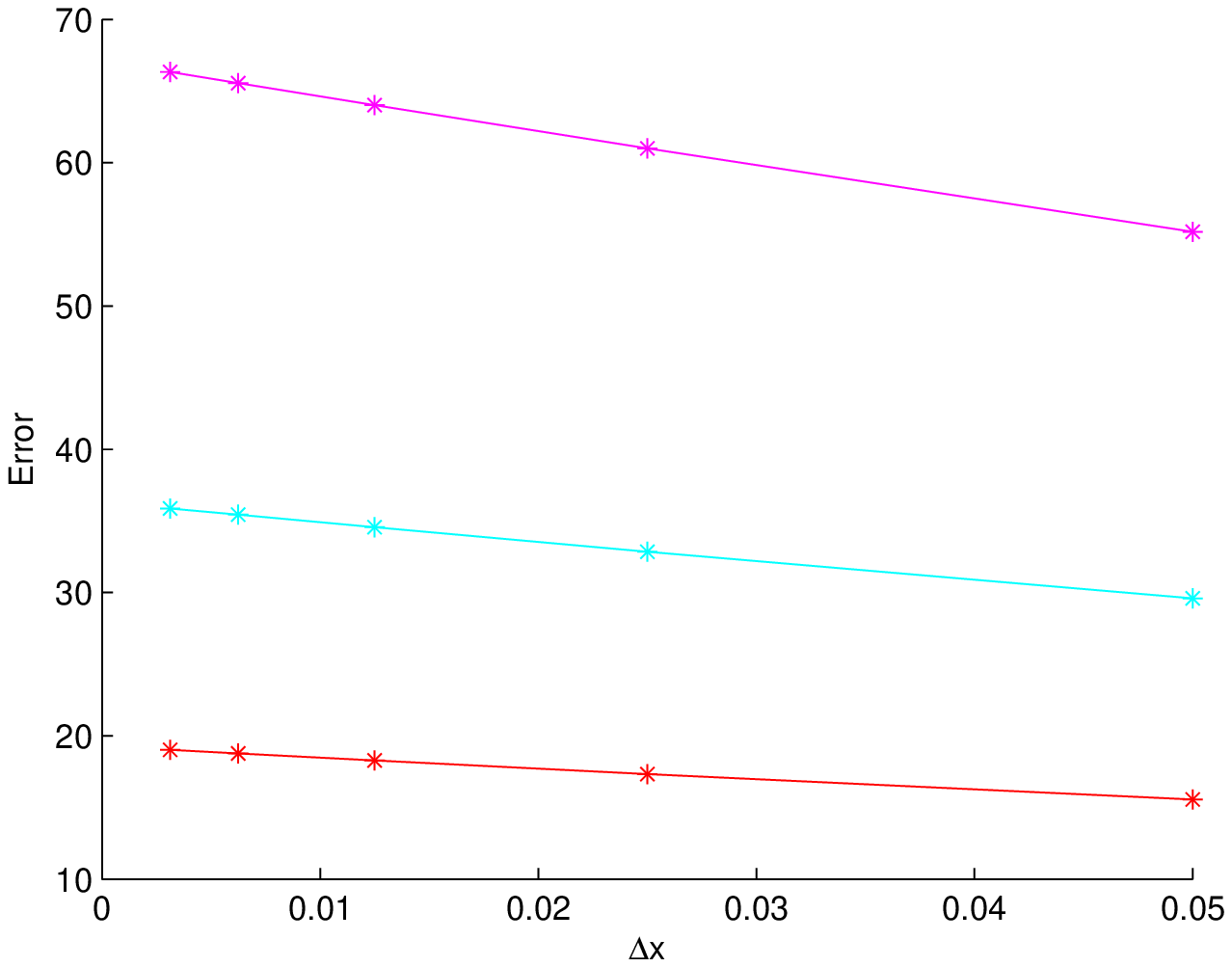}
\includegraphics[width=3.5cm,angle=-0]{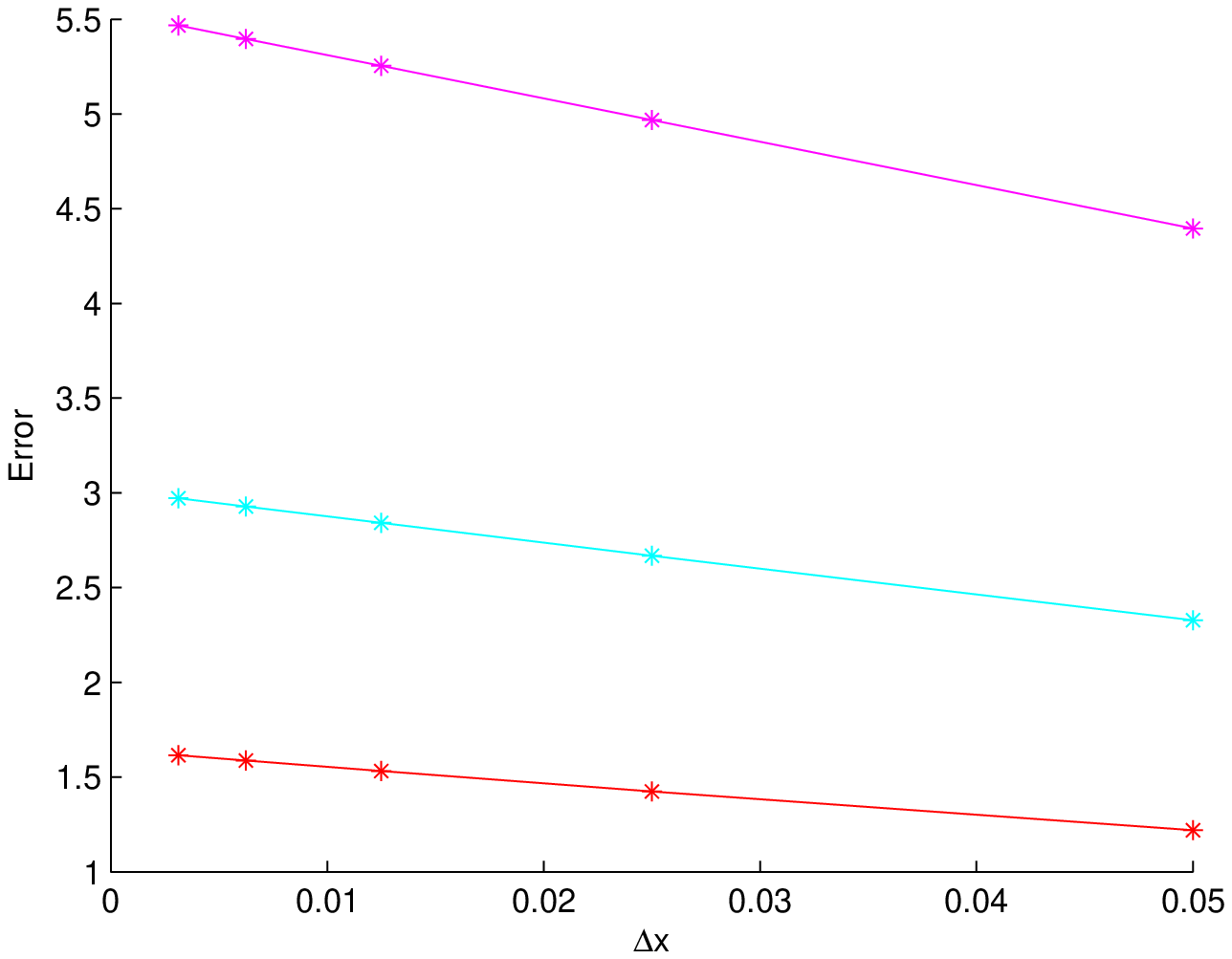} 
\includegraphics[width=3.5cm,angle=-0]{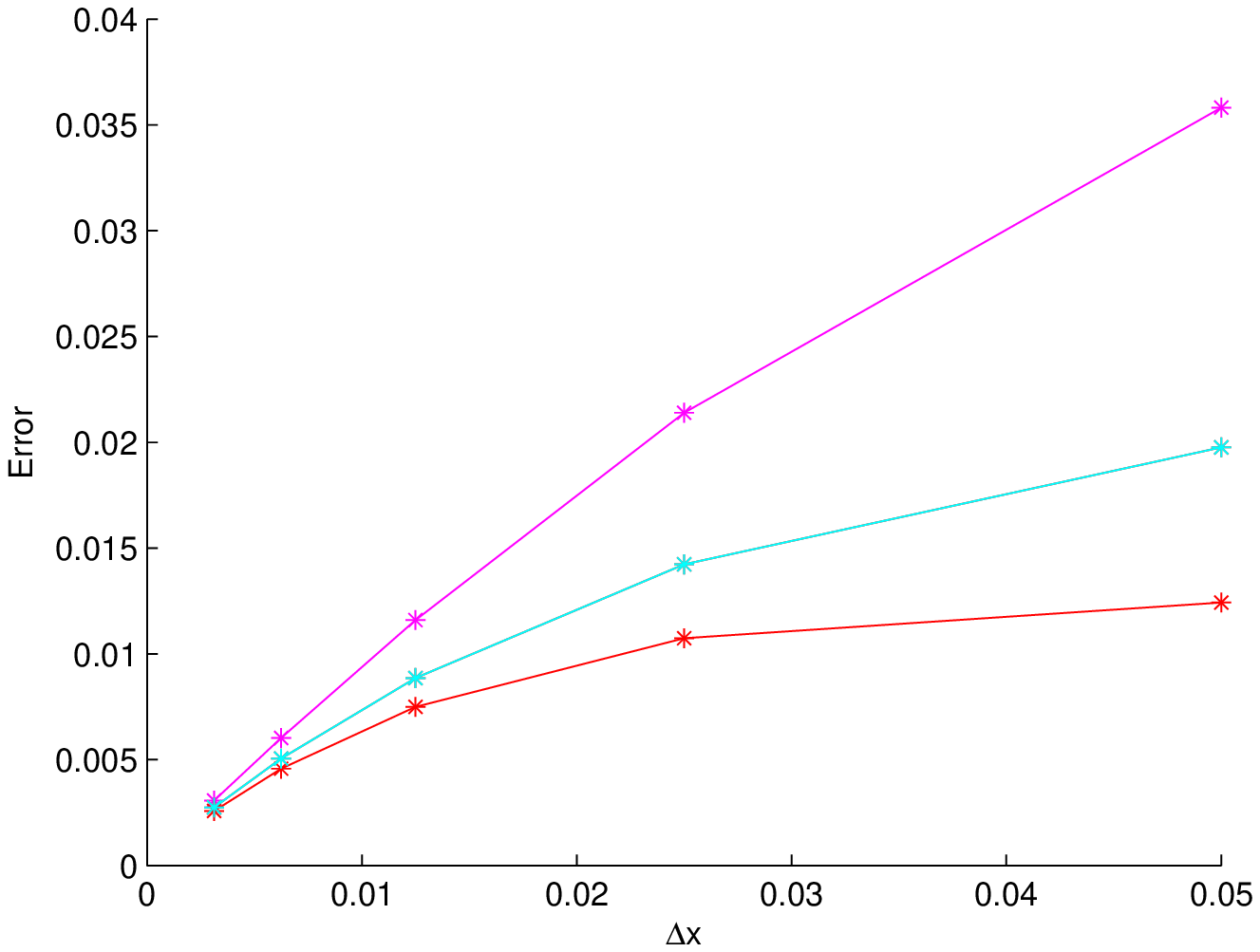} 
\end{center}
\caption{\label{figure_4} The solutions of the 1D Fisher's equation with the numerical splitting schemes for $K = 1, \; 0.5, \; \mbox{and} \; 0.25$ (asterisks in magenta $K=1.0$, cyan $K=0.5$ and red $K = 0.25$),  (left figure: $L_1$-norm, middle figure: $L_2$-norm and right figure: $L_{\infty}$-norm.}
\end{figure}

The results of the different norms are shown in Figure \ref{figure_4}.

\begin{remark}
We see the different solutions of the analytical and numerical solutions.
The numerical solutions fit to the analytical solutions. The nonlinearity is optimal resolved and more effective with the successive approximation method. 
Their resolution over the different time-dependent terms is more resolved as for Magnus-expansion, which only average their nonlinear properties.
\end{remark}

\subsection{Third Example: 2d Fisher's equation}

We deal with the 2D Fisher Equation, which describe the
spreading of genes see \cite{fisher1937} and has found applications in different
fields of research ranging from ecology \cite{kolomgorov1991} to plasma
physics \cite{gilding2004}.

We deal with a nonlinear PDE and split it into linear and nonlinear operators,
while we can compare to a analytical solution.

The Fisher's equation is given as
\begin{eqnarray}
\frac{\partial u(x, y, t)}{\partial t} & = & D_x \partial_{xx} u + D_y \partial_{yy} u + r (1 - \frac{u}{K}) u \; , \; \mbox{in} \; (x, t) \in \Omega \times [0, T] , \\
u(x, y, 0) &=& g(x, y) \; , \; \mbox{on} \; (x, y) \in \Omega , \\
u(x, y, t) &=& 0 \;   , \; \mbox{on} \; (x, y,  t) \in \partial \Omega \times [0, T],
\end{eqnarray}
%
where we assume $\Omega = [-2, 2] \times [-2, 2] \subset \R^2$,
such that we could apply the Dirichlet-boundary conditions 
$u(x, y, t) = 0 , \; \partial \Omega \times [0,T]$.
$u$ is the solution function, the initial condition is $g(x,y)$.
Form the dynamical view-point, we apply a homogeneous medium with $D$ as diffusion coefficient and we embed a growth of a logistic function, see \cite{vander2010}
with the $r$ is the growth rate and $K$ is the carrying capacity.

%
%

We apply the case $g(x, y) = \exp(-x^2 - y^2)$ .


We rewrite the equation-system (\ref{fisher_1}) in
operator notation, and obtain the following equations :
\begin{eqnarray}
\label{num_8}
  \partial_t u & = & A u + B(u) u \; , \\
   u(0) & = & u_0 ,
\end{eqnarray}
and we split our operators to a linear and nonlinear one:
\begin{eqnarray}
\label{num_9}
A =   D_x \partial_{xx} + D_y \partial_{yy} + r \\
B(u) = - r \frac{u}{K}  ,
\end{eqnarray}
with $D_x = D_y = 0.01, \; r = 1, K = 1$, later we apply the multiscale case with $K = 0.5, 0.25$.

%
%
%

%
%


Numerical Solution of the Diffusion part: \\

The operator $A$ is discretized as:
\begin{eqnarray}
A & = &  \frac{D}{\Delta x^2}\cdot  \left(\begin{array}{rrrrr}
 - T & I & ~ & ~ & ~ \\
 I & - T & I & ~ & ~ \\
 ~ & \ddots & \ddots & \ddots & ~ \\
 ~ & ~ & I & - T & I \\
 ~ & ~ & ~ & I & - T
\end{array}\right) \nonumber \\ 
& + & r  \left(\begin{array}{rrrrr}
 1 & 0 & ~ & ~ & ~ \\
 0 & 1 & ~ & ~ & ~ \\
 ~ & \ddots & \ddots & ~ & ~ \\
 ~ & ~ & 0 & 1 & 0 \\
 ~ & ~ & ~ & 0 & 1
\end{array}\right)~\in~\R^{I^2 \times I^2},
\end{eqnarray}
\begin{eqnarray}
T & = &  \left(\begin{array}{rrrrr}
4  & - 1 & ~ & ~ & ~ \\
-1 & 4 & -1 & ~ & ~ \\
 ~ & \ddots & \ddots & \ddots & ~ \\
 ~ & ~ & -1 & 4 & -1 \\
 ~ & ~ & ~ & -1 & 4
\end{array}\right) ~\in~\R^{I \times I},
\end{eqnarray}
where $I$ are the number of discretization points, e.g. $I = 100$.

The operator $B$ is written in the operator notation as:
\begin{eqnarray}
B(u) & = &  - \frac{r}{K}  \left(\begin{array}{rrrrr}
 u_1(t) & 0 & ~ & ~ & ~ \\
 0 & u_2(t) & ~ & ~ & ~ \\
 ~ & \ddots & \ddots & ~ & ~ \\
 ~ & ~ & 0 & u_{I^2-1}(t) & 0 \\
 ~ & ~ & ~ & 0 & u_{I^2}(t)
\end{array}\right)~\in~\R^{I^2 \times I^2},
\end{eqnarray}
where $I$ are the number of discretization points, e.g. $I = 100$.

The solution is given as ${\bf u}(t) = (u_1(t), \ldots, u_I(t), u_{I+1}(t), \ldots, u_{I^2}(t)$ . 
where $u_i(t) = u(x_j, y_k, t), x_j = \Delta x \; j, \; j = 1, \ldots, I, \; y_k = \Delta x \; j, \; k = 1, \ldots, I $ and $i = I (k-1) + j$.

Now, we can apply the discretized operator equations given in the following:
\begin{eqnarray}
\frac{\partial {\bf u}(t)}{\partial t} & = & A {\bf u} + B({\bf u}) {\bf u} \; , \; \mbox{in} \; t \in [0, T] \\
{\bf u}(0) &=& (g(x_1, y_1), \ldots, g(x_I, y_I))^t ,
\end{eqnarray}
 to our splitting schemes.


The nonlinear behaves in a non-commutative
manner $A' B - B' A \neq 0$.

We apply finite differences or finite elements to the spatial operator.
With $\Omega = [-2, 2] \times [-2, 2]$, $\Delta x = \Delta y = 0.05$ \\

We have the following results:

We have the following error in an $L_{\infty}$ Banach space:
\begin{eqnarray}
\label{kap7_gleich3}
E_{L_{\infty, \Delta x}}(t) & = &  \sup\nolimits_{x \in \Omega_h} | u_{num, \Delta x/8}(x, t) - u_{num}(x, t) |   \nonumber \\
              & = &   \sup\nolimits_{i = 1}^{N} | u_{num, \Delta x/8}(x_i, t) - u_{num}(x_i, t) |   .
\end{eqnarray}

We have the following error in an $L_1$ Banach space:
\begin{eqnarray}
\label{kap7_gleich3}
E_{L_{1, \Delta x}}(t) & = & \int_{\Omega_h} | u_{num, \Delta x/8}(x, t) - u_{num}(x, t) | \; dx  \nonumber \\
              & = & \Delta x \sum_{i=1}^N | u_{num, \Delta x/8}(x_i, t) - u_{num}(x_i, t) |  .
\end{eqnarray}

We have the following error in an $L_2$ Banach space:
\begin{eqnarray}
\label{kap7_gleich3}
E_{L_{2, \Delta x}}(t) & = & \sqrt{\int_{\Omega_h} (u_{num, \Delta x/8}(x, t) - u_{num}(x, t))^2 \; dx } \nonumber \\
              & = & \sqrt{\Delta x\sum_{i=1}^N (u_{num, \Delta x/8}(x_i, t) - u_{num}(x_i, t))^2  } .
\end{eqnarray}

\begin{figure}[ht]
\begin{center}  
\includegraphics[width=5.0cm,angle=-0]{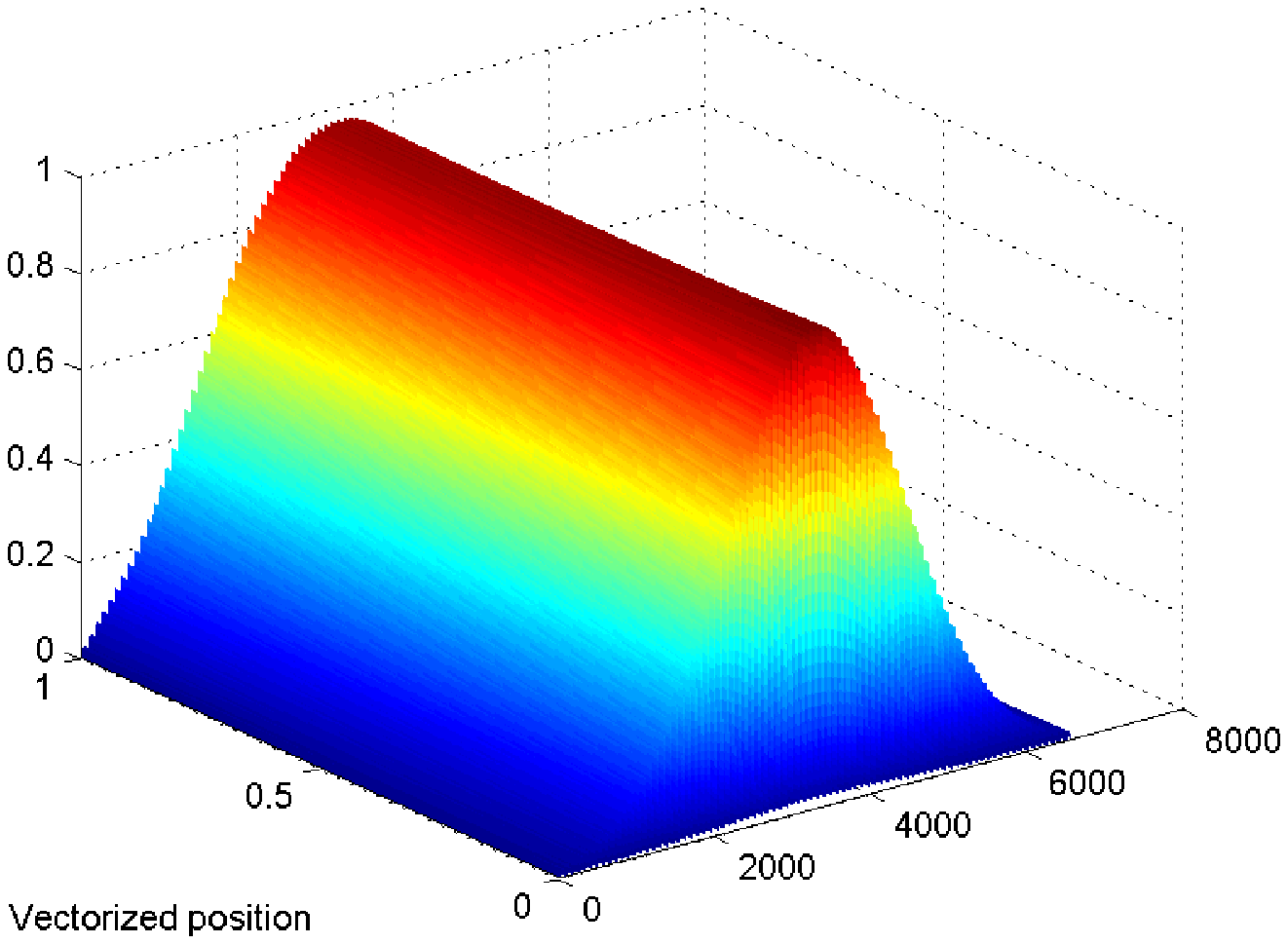} 
\end{center}
\caption{\label{figure_1_3} The solutions of the 2d Fisher's equation of the numerical solution.}
\end{figure}

The results of the different schemes are shown in Figure \ref{figure_1_3}.

\begin{figure}[ht]
\begin{center}  
\includegraphics[width=3.5cm,angle=-0]{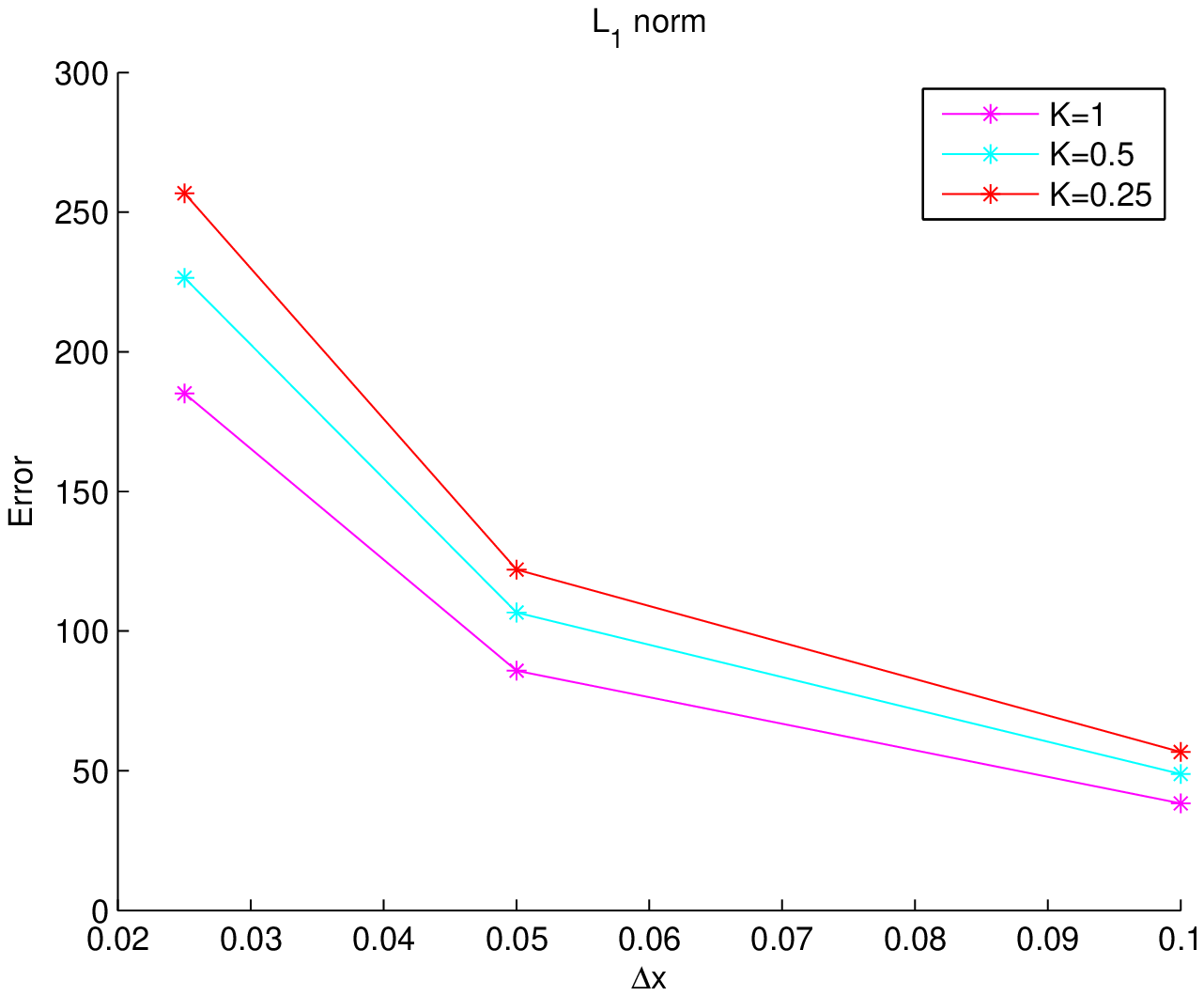} 
\includegraphics[width=3.5cm,angle=-0]{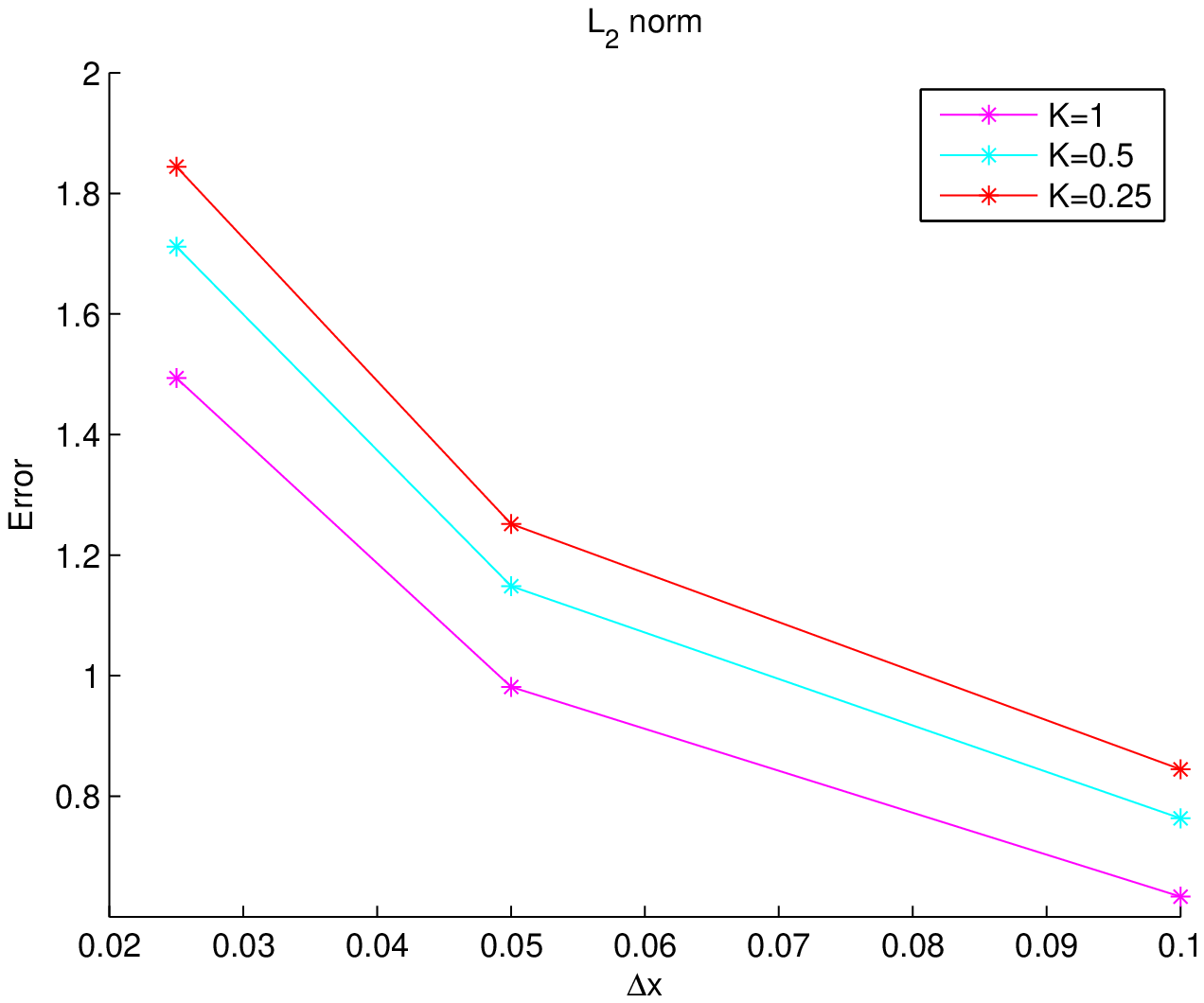} 
\includegraphics[width=3.5cm,angle=-0]{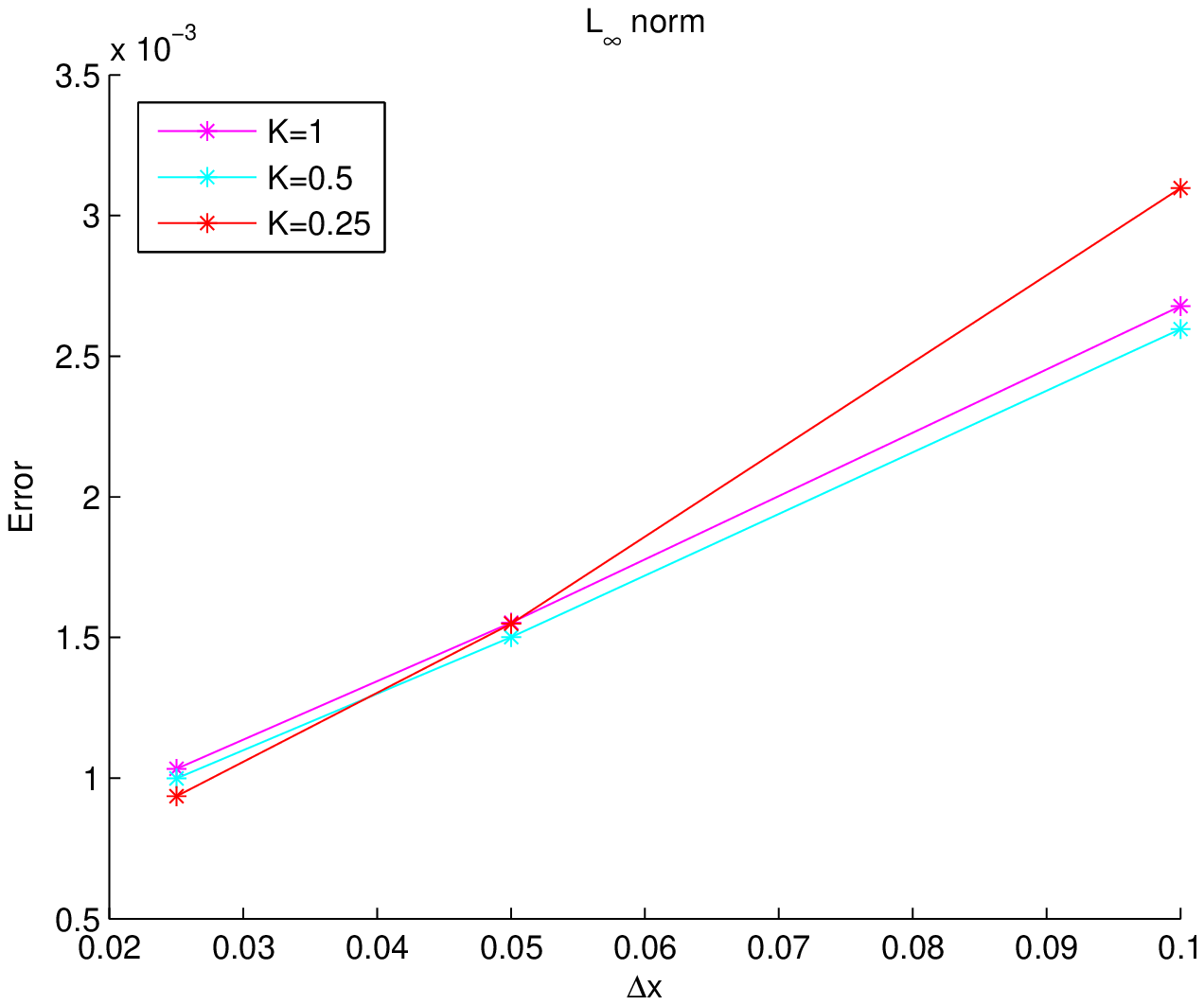} 
\end{center}
\caption{\label{figure_2_3} The convergence results of the 2d Fisher's equation (left figure: $L_1$ error, middle figure: $L_2$ error, right figure: $L_{inf}$ errors).}
\end{figure}

The convergence results of the different schemes are shown in Figure \ref{figure_2_3}.

\begin{remark}
We see the different solutions of the analytical and numerical solutions.
Here, we have a higher inverstigation to the spatial discretization and
therefore also for the solvers.
The numerical solutions fit to the analytical solutions. 
We have the same improvements also for the higher dimensions,
that we resolved the nonlinear terms via Taylor-expension in the
Multiscale method more accurate as with the Magnus-expansion.
\end{remark}

\subsection{Forth Example: 3d Fisher's equation}

We deal with the 3D Fisher Equation, which describe the
spreading of genes see \cite{fisher1937} and has found applications in different
fields of research ranging from ecology \cite{kolomgorov1991} to plasma
physics \cite{gilding2004}.

We deal with a nonlinear PDE and split it into linear and nonlinear operators,
while we can compare to a analytical solution.

The Fisher's equation is given as
\begin{eqnarray}
\frac{\partial u(x, y, z,  t)}{\partial t} & = & D_x \partial_{xx} u  + D_y \partial_{yy} u + D_z \partial_{zz} u \nonumber \\
&& + r (1 - \frac{u}{K}) u \; , \; \mbox{in} \; (x, t) \in \Omega \times [0, T] , \\
u(x, y, z, 0) &=& g(x, y, z) \; , \; \mbox{on} \; (x, y) \in \Omega , \\
u(x, y, z, t) &=& 0 \;   , \; \mbox{on} \; (x, y, z, t) \in \partial \Omega \times [0, T],
\end{eqnarray}
%
where we assume $\Omega = [-0.5, 0.5] \times [-0.5, 0.5] \times [-0.5, 0.5] \subset \R^3$,
such that we could apply the Dirichelt-boundary conditions 
$u(x, y, z, t) = 0 , \; \partial \Omega \times [0,T]$.
$u$ is the solution function, the initial condition is $g(x, y, z)$.
Form the dynamical view-point, we apply a homogeneous medium with $D$ as diffusion coefficient and we embed a growth of a logistic function, see \cite{vander2010}
with the $r$ is the growth rate and $K$ is the carrying capacity.

%
%

We apply the case $g(x, y) = \exp(-x^2 - y^2 - z^2)$ .


We rewrite the equation-system (\ref{fisher_1}) in
operator notation, and obtain the following equations :
\begin{eqnarray}
\label{num_8}
  \partial_t u & = & A u + B(u) u \; , \\
   u(0) & = & u_0 ,
\end{eqnarray}
and we split our operators to a linear and nonlinear one:
\begin{eqnarray}
\label{num_9}
A =   D_x \partial_{xx} + D_y \partial_{yy} + D_{zz} \partial_{zz} + r \\
B(u) = - r \frac{u}{K}  ,
\end{eqnarray}
with $D_x = D_y = D_z = 0.01, \; r = 1, K = 1$, later we apply the multiscale case with $K = 0.5, 0.25$.

%
%
%

%
%


Numerical Solution of the Diffusion part: \\

The operator $A$ is discretized as:
\begin{eqnarray}
A & = &  \frac{D}{\Delta x^2}\cdot  \left(\begin{array}{rrrrr}
 T_1 & I_1 & & \ldots  \\
 I_1 & T_1 & I_1 & ~ & ~ \\
 ~ & \ddots & \ddots & \ddots & ~ \\
 ~ & ~ & I_1 & T_1 & I_1 \\
 ~ & ~ & ~ & I_1 & T_1
\end{array}\right)  \nonumber \\
 & + & r  \left(\begin{array}{rrrrr}
 1 & 0 & ~ & ~ & ~ \\
 0 & 1 & ~ & ~ & ~ \\
 ~ & \ddots & \ddots & ~ & ~ \\
 ~ & ~ & 0 & 1 & 0 \\
 ~ & ~ & ~ & 0 & 1
\end{array}\right)~\in~\R^{I^3 \times I^3},
\end{eqnarray}
\begin{eqnarray}
T_1 & = &  \left(\begin{array}{rrrrr}
  T & I & & \ldots  \\
 I & T & I & ~ & ~ \\
 ~ & \ddots & \ddots & \ddots & ~ \\
 ~ & ~ & I & T & I \\
 ~ & ~ & ~ & I & T
\end{array}\right) ~\in~\R^{I^2 \times I^2},
\end{eqnarray}
\begin{eqnarray}
I_1 & = &  \left(\begin{array}{rrrrr}
  I & 0 & & \ldots  \\
 0 & I & 0 & ~ & ~ \\
 ~ & \ddots & \ddots & \ddots & ~ \\
 ~ & ~ & = & I & 0 \\
 ~ & ~ & ~ & 0 & I
\end{array}\right) ~\in~\R^{I^2 \times I^2},
\end{eqnarray}

\begin{eqnarray}
T & = &  \left(\begin{array}{rrrrr}
-6  & 1 & ~ & ~ & ~ \\
1 & -6 & 1 & ~ & ~ \\
 ~ & \ddots & \ddots & \ddots & ~ \\
 ~ & ~ & 1 & -6 & 1 \\
 ~ & ~ & ~ & 1 & -6
\end{array}\right) ~\in~\R^{I \times I},
\end{eqnarray}

\begin{eqnarray}
I & = &  \left(\begin{array}{rrrrr}
1  & 0 & ~ & ~ & ~ \\
0 & 1 & 0 & ~ & ~ \\
 ~ & \ddots & \ddots & \ddots & ~ \\
 ~ & ~ & 0 & 1 & 0 \\
 ~ & ~ & ~ & 0 & 1
\end{array}\right) ~\in~\R^{I \times I},
\end{eqnarray}
where $I$ are the number of discretization points, e.g. $I = 100$.

The operator $B$ is written in the operator notation as:
\begin{eqnarray}
B(u) =  - \frac{r}{K}  \left(\begin{array}{rrrrr}
 u_1(t) & 0 & ~ & ~ & ~ \\
 0 & u_2(t) & ~ & ~ & ~ \\
 ~ & \ddots & \ddots & ~ & ~ \\
 ~ & ~ & 0 & u_{I^3-1}(t) & 0 \\
 ~ & ~ & ~ & 0 & u_{I^3}(t)
\end{array}\right)~\in~\R^{I^3 \times I^3},
\end{eqnarray}
where $I$ are the number of discretization points, e.g. $I = 100$.

The solution is given as ${\bf u}(t) = (u_1(t), \ldots, u_I(t), u_{I+1}(t), \ldots, u_{I^3}(t)$ . 
where $u_i(t) = u(x_j, y_k, z_l, t), x_j = \Delta x \; j, \; j = 1, \ldots, I, \; y_k = \Delta x \; j, \; k = 1, \ldots, I \; y_l = \Delta x \; l, \; l = 1, \ldots, I $ and $i = I^2 (l-1) + I (k-1) + j$.

Now, we can apply the discretized operator equations to our splitting schemes,
which are given as:
\begin{eqnarray}
\frac{\partial {\bf u}(t)}{\partial t} & = & A {\bf u} + B({\bf u}) {\bf u} \; , \; \mbox{in} \; t \in [0, T] \\
{\bf u}(0) &=& (g(x_1, y_1, z_1), \ldots, g(x_I, y_I, z_I))^t .
\end{eqnarray}
%


The nonlinear behaves in a non-commutative
manner $A' B - B' A \neq 0$.

We apply finite differences or finite elements to the spatial operator.
With  $\Omega = [-0.5, 0.5] \times [-0.5, 0.5] \times [-0.5, 0.5], \Delta x = \Delta y = \Delta z = 0.05$.

We have the following results:

We have the following relative error in an $L_{\infty}$ Banach space:
\begin{eqnarray}
\label{kap7_gleich3}
E_{L_{\infty, \Delta x},rel}(t) & = & \frac{\sup\nolimits_{x \in \Omega_h} | u_{num, \Delta x/8}(x, t) - u_{num}(x, t) |}{\sup\nolimits_{x \in \Omega_h} | u_{num, \Delta x/8}(x, t) |}   \nonumber \\
              & = &  \frac{ \sup\nolimits_{i = 1}^{N} | u_{num, \Delta x/8}(x_i, t) - u_{num}(x_i, t) | } {\sup\nolimits_{i = 1}^{N} | u_{num, \Delta x/8}(x_i, t) |}   .
\end{eqnarray}

We have the following error in an $L_1$ Banach space:
\begin{eqnarray}
\label{kap7_gleich3}
E_{L_{1, \Delta x}, rel}(t) & = & \frac{\int_{\Omega_h} | u_{num, \Delta x/8}(x, t) - u_{num}(x, t) | \; dx}{\int_{\Omega_h} | u_{num, \Delta x/8}(x, t)| \; dx}  \nonumber \\
              & = & \frac{\Delta x \sum_{i=1}^N | u_{num, \Delta x/8}(x_i, t) - u_{num}(x_i, t) |}{\Delta x \sum_{i=1}^N | u_{num, \Delta x/8}(x_i, t)|}  .
\end{eqnarray}

We have the following error in an $L_2$ Banach space:
\begin{eqnarray}
\label{kap7_gleich3}
E_{L_{2, \Delta x,rel}}(t) & = & \frac{\sqrt{\int_{\Omega_h} (u_{num, \Delta x/8}(x, t) - u_{num}(x, t))^2 \; dx }} { \sqrt{\int_{\Omega_h} (u_{num, \Delta x/8}(x, t))^2 \; dx } } \nonumber \\
              & = & \frac{\sqrt{\Delta x\sum_{i=1}^N (u_{num, \Delta x/8}(x_i, t) - u_{num}(x_i, t))^2  } }{\sqrt{\Delta x\sum_{i=1}^N (u_{num, \Delta x/8}(x_i, t) } } .
\end{eqnarray}

\begin{figure}[ht]
\begin{center}  
\includegraphics[width=3.5cm,angle=-0]{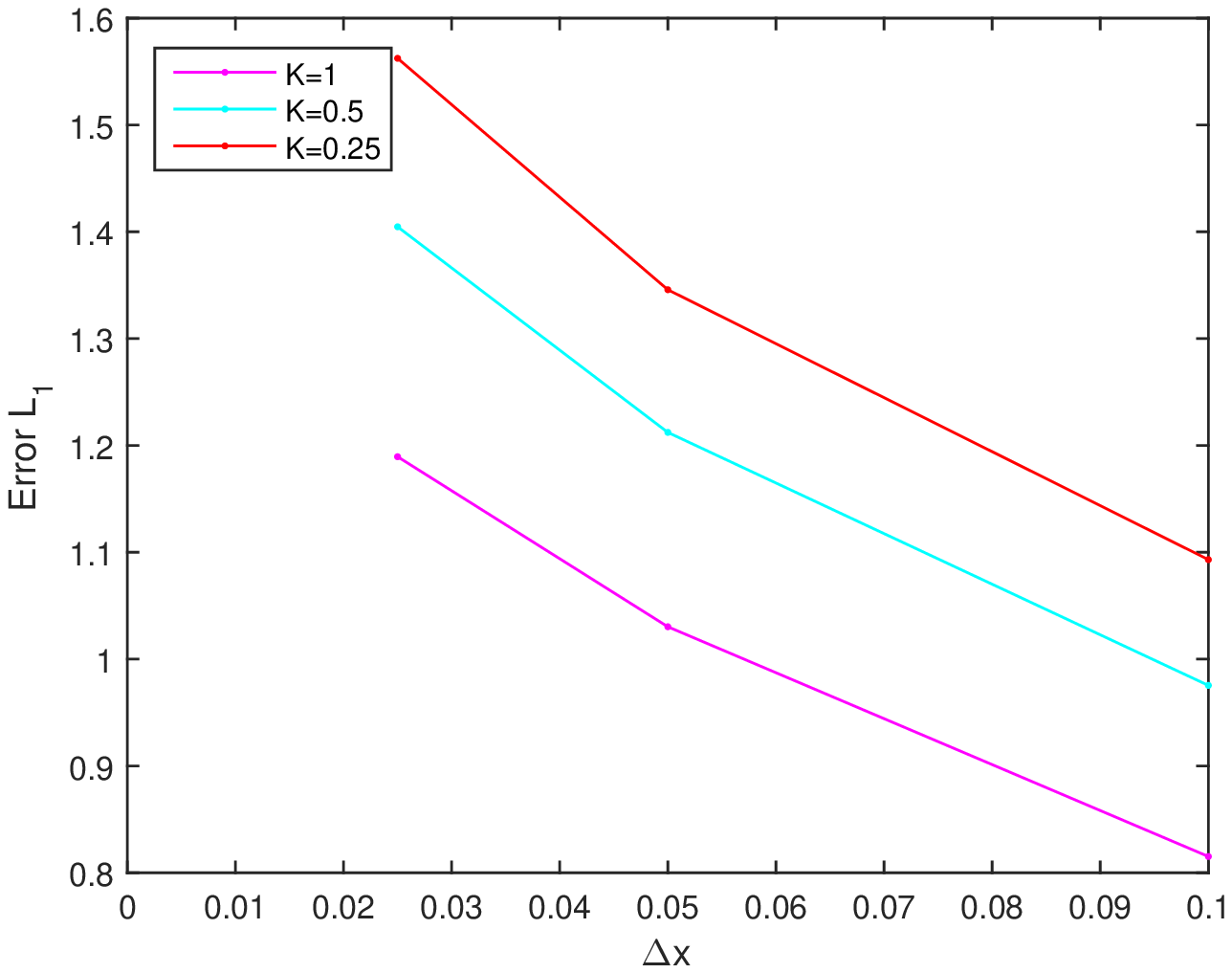} 
\includegraphics[width=3.5cm,angle=-0]{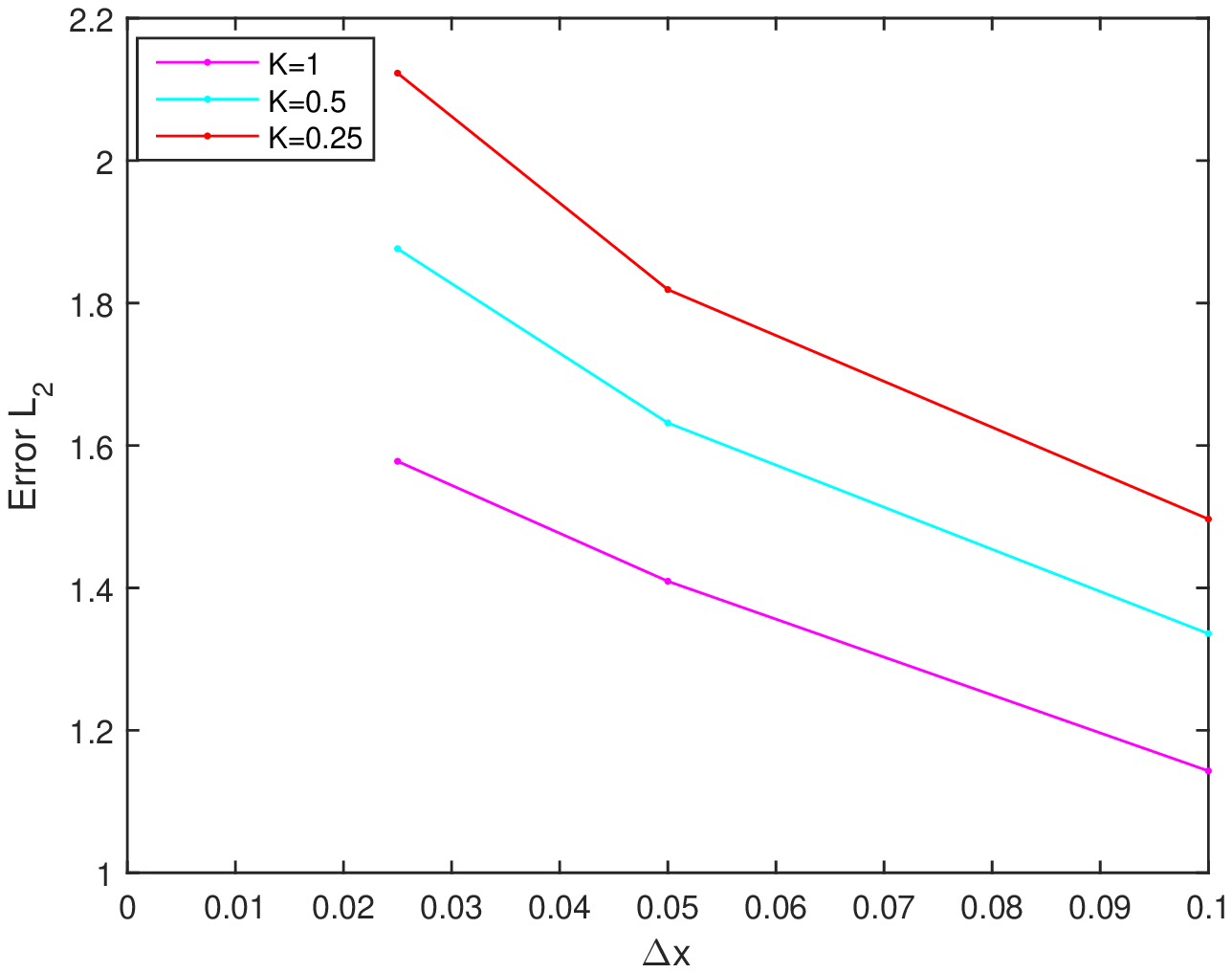} 
\includegraphics[width=3.5cm,angle=-0]{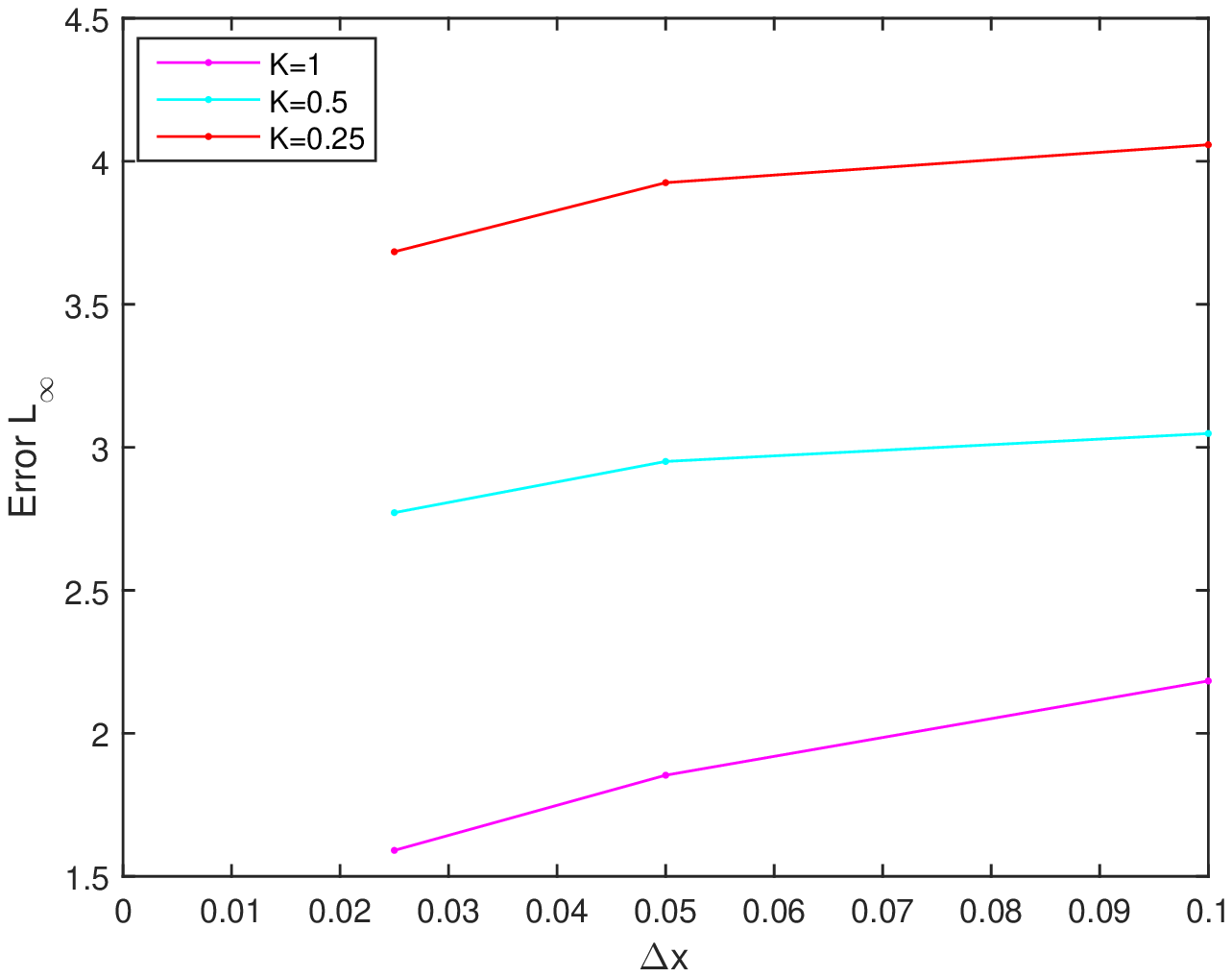} 
\end{center}
\caption{\label{figure_3_4} The convergence results of the 3d Fisher's equation (left figure: $L_1$ relative error, middle figure: $L_2$ relative error, right figure: $L_{inf}$ relative errors).}
\end{figure}

The convergence results of the different schemes are shown in Figure \ref{figure_3_4}.

\begin{remark}
We see the different solutions of the analytical and numerical solutions.
Here, we have a higher investigation to the spatial discretization and
apply a fast computation via Leja-point of the $exp$-matrices,
see \cite{caliari2014}.
The numerical solutions fit to the analytical solutions. 
We have the same improvements also for the higher dimensions, here we
also resolve the nonlinear terms  more accurate as with the Magnus-expansion.
\end{remark}

\section{Conclusions and Discussions }
\label{concl}

In this work, we have presented application to successive approximations
that are related to iterative splitting schemes.
We present the convergence analysis of the scheme and approximation
to multiple scale methods. We see the benefits in resolving the nonlinearity
in the Taylor-expension, such that we could conclude a higher order approximation of the nonlinearity in the recent time-step. Numerical experiments present the benefit of the scheme to standard Magnus expension methods. In future, we analyse and apply our method to different real-life applications.

\end{document}